\documentclass[12pt]{amsart}
\usepackage{amsbsy}
\usepackage{graphicx,epsfig,subfigure,psfrag}
\begin{document}

\newtheorem{theorem}{Theorem}
\newtheorem{proposition}{Proposition}
\newtheorem{lemma}{Lemma}
\newtheorem{corollary}{Corollary}
\newtheorem{definition}{Definition}
\newtheorem{remark}{Remark}
\newcommand{\beq}{\begin{equation}}
\newcommand{\eeq}{\end{equation}}
\numberwithin{equation}{section}
\numberwithin{theorem}{section}
\numberwithin{proposition}{section}
\numberwithin{lemma}{section}
\numberwithin{corollary}{section}
\numberwithin{definition}{section}
\numberwithin{remark}{section}
\newcommand{\re}{{\mathbb R}}
\newcommand{\n}{\nabla}
\newcommand{\ren}{{\mathbb R}^N}
\newcommand{\iy}{\infty}
\newcommand{\pa}{\partial}
\newcommand{\fp}{\noindent}
\newcommand{\ms}{\medskip\vskip-.1cm}
\newcommand{\mpb}{\medskip}
\newcommand{\BB}{{\bf B}}
\newcommand{\Am}{{\bf A}_{2m}}
\newcommand{\ssk}{\smallskip}
\renewcommand{\a}{\alpha}
\renewcommand{\b}{\beta}
\newcommand{\g}{\gamma}
\newcommand{\G}{\Gamma}
\renewcommand{\d}{\delta}
\newcommand{\D}{\Delta}
\newcommand{\e}{\varepsilon}
\renewcommand{\l}{\lambda}
\renewcommand{\o}{\omega}
\renewcommand{\O}{\Omega}
\newcommand{\s}{\sigma}
\renewcommand{\t}{\tau}
\renewcommand{\th}{\theta}
\newcommand{\z}{\zeta}
\newcommand{\wx}{\widetilde x}
\newcommand{\wt}{\widetilde t}
\newcommand{\noi}{\noindent}
\newcommand{\supp}{{\rm supp}\,}
\def\com#1{\fbox{\parbox{6in}{\texttt{#1}}}}

\title 
{\bf On self-similar collapse of discontinues data
 for thin film equations
 with doubly degenerate mobility}

\author{V.A.~Galaktionov}

\address{Department of Mathematical Sciences, University of Bath,
 Bath BA2 7AY, UK}

\email{vag@maths.bath.ac.uk}


 \keywords{Thin film equation, degenerate mobility, Riemann's problem, similarity solution,
 the Cauchy problem, interface equations, oscillatory solutions.
 }
 \subjclass{35K55, 35K40, 35K65. }
\date{\today}




\begin{abstract}
As a basic model,  the fourth-order quasilinear thin film equation
(the TFE--4) with a concentration-dependent doubly degenerate
mobility coefficient vanishing  at two equilibrium levels $u = \pm
1$,
 $$
u_t =- (|1-u^2| u_{xxx})_x \, ,
 $$
 is
 studied.
The basic {\em Riemann problem} for the TFE--4
 with discontinuous initial data
  $$
 u(x,0)= {\rm sign} \, x= \left\{ \begin{matrix}\,\,\,  1, \,\,\, x>0, \\
 -1, \,\,\, x<0,
 \end{matrix}
 \right.
 $$
  is considered.
This problem is shown to  admit self-similar solutions of the form
 $$
  \begin{matrix}
 u_+(x,t)= f(y), \quad \mbox{where} \quad y=x/t^{\frac 14}
 \,\,\, \mbox{and $f$ solves the ODE}\smallskip\smallskip\\
   -(|1-f^2| f''')' + \frac 14 \, f' y=0, \quad f(\pm \infty)= \pm
   1.
   \end{matrix}
 $$
  The similarity
 profiles are  different for the zero contact angle and zero-flux  FBP in
  a bounded domain and for
 the Cauchy problem in $\re \times \re_+$.
  The similarity profiles can be also
  obtained by branching at $n=0$ for the mobility coefficient
  $|1-f^2|^n$ that involves the linear bi-harmonic equation
   $
   u_t=-u_{xxxx}.
   $

Finite propagation in the same Riemann problem is studied for the
TFE--4 with the unstable diffusion term such as
 $$
 u_t =-
(|1-u^2| u_{xxx})_x -
 (|1-u^2|^m u_x)_{x},
 $$
  where $m \ge 0$ is a parameter,
    and dynamic interface
 equations are derived.
  Similar conclusions on similarity solutions of Riemann's problem
  apply
   to the $2m$th-order TFEs,
 $$
u_t =(-1)^{m+1} (|1-u^2|^n D_x^{2m-1} u)_x \quad (m \ge 2 \,\,
\mbox{integer}, \,\,\, n>0).
 $$

\end{abstract}

\maketitle


\setcounter{equation}{0}
\section{Introduction: Riemann's problems
 and
similarity solutions}
 \label{Sect1}
  \setcounter{equation}{0}

We study the collapse of initial discontinuities (singularities)
in higher-order quasilinear degenerate parabolic equations of thin
film type. Such singularity problems are classical in PDE theory
and systematically  occurred in entropy mathematics of
conservation laws and hyperbolic
 systems in the 1950s.

\subsection{On Riemann-type problems in nonlinear evolution PDEs}

Riemann's problems play the truly special role  and are key for
many nonlinear evolution equations with partial derivatives, to be
discussed in
a short survey below. Typical Riemann's problems have been studied
for nonlinear PDEs with quasilinear spatial differential operators
of the  orders up to three that are denoted below as equation
classes
 {\bf (1)}, {\bf (2)}, and {\bf (3)}.

\smallskip

\noi\underline{\sc {\bf (1)} Conservation laws}. This is the most
classic matter associated with gas dynamics problems, where the
shock and blow-up Guderley-type regimes were actively studied
since the 1930s\footnote{First results on the formation of shocks
were due to Riemann himself  in 1858 \cite{Ri58}; see
\cite{Chr07}.}.
 Namely,
 {\em Riemann's problems} are  well-known  for scalar conservation
laws such as the first-order 1D {\em Euler equation} as the key
representative,
 \beq
 \label{E1}
 u_t + u u_x=0 \quad \mbox{in}
 \quad  \re \times \re_+.
  \eeq
  According to classic theory (see details in Smoller
  \cite[Part~III]{Sm}),
  for such first (and, {\em odd}-) -order PDEs, two Riemann's problems are key with
 discontinuous initial data
 \beq
 \label{ss11}
 S_-(x)= -{\rm sign} \, x \quad \mbox{and}\quad S_+(x)= {\rm sign} \,
 x.
  \eeq
The first initial function leads to the {\em shock wave}, which is
the entropy stationary solution $u_-(x,t) \equiv S_-(x)$ for all
$t> 0$, while the second one gives
 rise to the {\em rarefaction
wave}
 \beq
 \label{rr1}
   \mbox{$
   u_{\rm +}(x,t)= f\big( \frac x t \big) \equiv \left\{ \begin{matrix}
    -1 \quad \mbox{for} \,\,\, x<-t, \cr
 \,\, \frac x t \quad  \mbox{for} \,\,\, |x|<t, \cr
  1 \quad \,  \mbox{for} \,\,\, x>t.
 \end{matrix}
   \right.
    $}
 \eeq
The above formula (\ref{rr1}) shows the generic  {\em similarity}
mechanism, with the scaling invariant $\frac x t$, of the collapse
of this non-entropy initial singularity $S_+(x)$. For the first
data in (\ref{ss11}), replacing $t \mapsto -t$, $f \mapsto -f$ in
(\ref{rr1}), gives the formation of the shock $S_-(x)$ in a
blow-up singularity as $t \to 0^-$ via the {\em self-similar
solution}, with the same profile $f$,
 \beq
 \label{rr2}
  \mbox{$
 u_-(x,t)=- f \big(\frac x{(-t)}\big) \to S_-(x) \quad \mbox{as} \quad t \to 0^-.
 $}
  \eeq
Riemann's problems describing   shock wave formation and collapse
of other discontinuities are also classic for strictly hyperbolic
systems since the 1950s; see Bressan \cite{Bres} and Dafermos
\cite{Daf} for history and main results.

\smallskip

\noi\underline{\sc {\bf (2)} Quasilinear heat equations}. The
collapse of initial discontinuity $S_+(x)$ for the second-order
parabolic equations such as
 \beq
 \label{rr3}
 u_t=(k(u)u_x)_x,
  \eeq
  where $k(u) \ge 0$ is a given coefficient ($k(u)=|u|$  corresponds to  the {\em signed  porous
  medium equation}), is also described by {\em similarity} patterns
  with the scaling invariant $\frac x{\sqrt t}$ (sometimes called the Bolzman
  substitution\footnote{Similarity
solutions were used by Weierstrass around 1870, and by Bolzman
around 1890. In parabolic PDEs, this rescaled variable $y$,
together with the backward blow-up one $y=\frac x{\sqrt{-t}}$ for
$t<0$, was freely and effectively used by Sturm in 1836
\cite{St36} in his classification of all possible multiple zeros
of $u(x,t)$.})
   \beq
   \label{rr4}
     \begin{matrix}
   u_+(x,t)= f(y), \quad y= \frac x{\sqrt t},
   \quad \mbox{where $f$ solves the ODE} \smallskip\smallskip\smallskip\\
    (k(f)f')'+ \frac 12\, f' y=0 \quad \mbox{in} \quad \re \quad (f(\pm \infty)= \pm
   1).\qquad\,\,
   \end{matrix}
    \eeq
    To our knowledge,
first results on existence of such similarity solutions were due
to Polubarinova-Kochina (1948) \cite{Pol}. Later, such ODEs have
been studied in detail for arbitrary nonlinearities $k(u)$; see
Atkinson--Peletier \cite{AP0, AP}.
 By parabolic interior  regularity, the opposite pattern $u_-(x,t)$,
  describing  formation of a shock, is
 nonexistent, i.e., discontinuous solutions cannot appear evolutionary in parabolic PDEs (\ref{rr3}).

\smallskip

\noi\underline{\sc {\bf (3)} Rarefaction and shock  waves for
NDEs}. Riemann-type problems on discontinuous solutions $S_-(x)$
and collapse of $S_+(x)$ are actual for third-order nonlinear
dispersion equations (NDEs--3), e.g.,
 such as
   \beq
 \label{ent.12}
 u_t=(uu_x)_{xx},
  \eeq
 which, besides its applications in various problems of  shallow
 water theory and compact pattern formation (see references in \cite[Ch.~4]{GSVR}), are a natural
 extension of the first-order conservation law (\ref{E1}) to
the third- and other odd-order cases. Actually, writing
(\ref{ent.12}) in a ``pseudo-differential" form
 \beq
 \label{ee1}
 P \, u_t + uu_x=0,
  \eeq
  with a standard definition of the inverse Laplacian
  $P=(-D_x^2)^{-1}>0$ (in a bounded or unbounded interval),
   (\ref{ent.12}) can be viewed as a {\em non-local
 conservation law}. Note that, by no means, (\ref{ent.12}) can be
 treated as a hyperbolic system.

 Concerning the rarefaction
wave for (\ref{ent.12}) with data $S_+(x)$,
 this
  is again
  described by the {\em similarity solution}
  \beq
  \label{ent.14}
   \mbox{$
  u_{+}(x,t) = f(y), \,\,\, \mbox{where} \,\, y=
  \frac x{t^{1/3}}
   $}
   \eeq
 and $f$ solves the following ODE obtained on substitution into (\ref{ent.12}):
  \beq
  \label{ent.15}
  \mbox{$
   (f f')''+ \frac 13 \, f'y=0
    \quad  (f(\pm \infty)= \pm 1).
     $}
     \eeq
 This is a more difficult ODE, though  general understanding and some
  results on existence of rarefaction
 (and also shock) profiles are available; see \cite[\S~7]{GalEng},
 \cite[p.~173]{GSVR},
  and \cite{GalNDE5, GPnde} for more advanced $\d$-entropy theory for the NDEs such as (\ref{ent.12})
  and higher-order ones.

 Similarly, finite-time blow-up formation as $ t \to 0^-$ of the shock $S_-(x)$ is
 given by the similarity solution
\beq
  \label{ent.14n}
   \mbox{$
  u_{-}(x,t) = f(y), \,\,\, \mbox{where} \,\, y=
  x/{(-t)^{1/3}},
   $}
   \eeq
 with the reflected profile
  $-f(y)$  satisfying the ODE  (\ref{ent.15}).


\subsection{(4) Riemann's problems for thin film models with degenerate mobility.}


Thus, as the next extension of the three classes of PDEs
characterized above, we naturally arrive at  the fourth-order
quasilinear parabolic PDEs, for which Riemann's problems, though
not leading to shock waves, become more difficult.

 In this paper,
our main model is the fourth-order {\em thin film equation} (the
TFE--4)
 with concentration-dependent doubly degenerate mobility
coefficient $b(u)=1-u^2$; see   \cite{CENC,EG96,BBG99} for further
applied motivations and various results. For convenience, in the
dimensionless form, we write it as follows:
 \beq
 \label{1}
 u_t =- (|1-u^2| u_{xxx})_x,
  \eeq
 so it  is degenerate at two levels $u = \pm 1$, along which
  one can expect finite propagation of interfaces.
 We pose  (\ref{1}) either in a bounded domain (the case of the zero
 contact angle and zero-flux FBP), or in the whole space $\re \times \re_+$
 (this we call the Cauchy problem,
  CP). Concerning further physical derivations, applied motivations, history, and discussion
  of related models, we refer
  to Novick-Cohen's most recent monograph
\cite{CHbook}.




The absolute value in the nonlinear coefficient  $|1-u^2|$ in
(\ref{1}) is introduced to cover the case of the Cauchy problem or
other FBPs, where, as we show, for initial data within the range,
$u_0(x) \in [-1,1]$, the solution becomes oscillatory  about $u =
\pm 1$ and gets out of the physical range $[-1,1]$. As usual in
TFE problems, this may cause difficulties in the physical
motivation of such a situation. Nevertheless, such solutions are
necessary for general TFE theory and for understanding the
difference between the Cauchy and free-boundary problems for
degenerate TFEs. In this case,   we are obliged  to take the
absolute value
 $|1-u^2|$ of the mobility coefficients in (\ref{1}) to keep
the equation degenerate {\em parabolic} for all $u \in \re$.
 Otherwise, the parabolic problem becomes partially backward in time
 and loses its classic local well-posedness. Such  instabilities of the mathematical nature
 can be  suspicious
 from the physical point of view, though sometimes should be taken
 into account and lead to Young's measure solutions (not properly developed for such
 higher-order TFEs).



 For the  zero-flux and zero
contact angle FBP and solutions $u \in [-1,1]$ from the range, one
can use the standard TFE--4
  \beq
   \label{eq1}
   u_t =- ((1-u^2) u_{xxx})_x.
   \eeq
Existence and some regularity results for weak solutions of
(\ref{1}) and similar PDEs with extra low-order terms are
well-known since the pioneering paper by Bernis and Friedman
\cite{BF1}; see also \cite{EG96} and results and references in
more recent papers
  \cite{BerPugh1, BerPugh2, Bl4, SPugh, WitBerBer}.

 \subsection{Two basic Riemann's problems for the TFE--4 and similarity
 reduction}

 In view of reflection symmetry $x \mapsto -x$ of (\ref{1}),
there exists the single Riemann problem for the TFE--4 (\ref{1})
with initial data
  \beq
  \label{2}
 u(x,0)= S_+(x)= {\rm sign} \, x= \left\{ \begin{matrix}\,\,\, 1, \,\,\, x>0, \\
 -1, \,\,\, x<0.
 \end{matrix}
 \right.
 \eeq
 We call this the Riemann-1 problem or simply the RP--1.


 By parabolic interior regularity, one can expect to have a sufficiently smooth solution
 of the FBP (\ref{eq1}), (\ref{2}) for $t>0$  as a
 kind of rarefaction wave in parabolic setting.
The collapse of such initial singularity and propagation of
interfaces for small $t>0$ are key for understanding multi-phase
phenomena within of framework of various nonlinear problems of
TFE, Allen--Cahn, and Cahn--Hilliard type.


For the sake of mathematical consistency (and not only that), we
also introduce another Riemann-type problem with a half of the
step in (\ref{2})  setting,
  \beq
  \label{2NN}
 u(x,0)=  \theta(-x) \equiv \left\{ \begin{matrix} 0, \,\,\,\, x>0, \\
 1, \,\,\, x<0,
 \end{matrix}
 \right.
 \eeq
 where $\th(x)$ is
 the Heaviside function.
We refer to this as the RP--$\frac 12$. Notice the principle
difference in comparison with the RP--1 in (\ref{2}): for
(\ref{2NN}) we expect finite propagation along the equilibrium
level $u=1$ only, while on $\{u=0\}$ the propagation is infinite.




\smallskip

Our main goal is to demonstrate that, in general, both Riemann's
problems (\ref{1}), (\ref{2}) and (\ref{2NN}) admit  solutions
with the self-similar structure
 \beq
 \label{3}
  \mbox{$
 u_+(x,t)= f(y), \quad \mbox{where} \quad y=\frac x{t^{ 1/4}}.
  $}
  \eeq
  Then for the RP--1,
  $f$ solves the boundary-value problem for the ODE
  which, by anti-symmetry, we pose in $\re_+$ only,
  \beq
  \label{4}
    \textstyle{
   {\bf B}(f) \equiv -(|1-f^2| f''')' + \frac 14 \, f' y=0 \quad \mbox{for} \quad
   y>0,
    }
  \eeq
 with two anti-symmetry conditions at the origin,
 \beq
 \label{5}
 f(0)=f''(0)=0.
 \eeq
More precisely, for the FBP setting, we need to find a {\em
finite} $y_0>0$ such that the zero contact angle and zero-flux
conditions hold,
 \beq
 \label{6}
  f(y_0)=1, \quad f'(y_0)= (|1-f^2| f''')(y_0)=0 \quad ({\rm FBP}).
   \eeq
Similarly, for the RP--$\frac 12$, with data (\ref{2NN}) we obtain
the ODE (\ref{3}) in $\re$ with the regular condition at
$y=-\infty$
 \beq
 \label{5R}
 f(-\infty)=0,
  \eeq
and the free-boundary conditions (\ref{6}).

For the CP, we consider equation (\ref{4})  with the single
condition
 \beq
 \label{7}
 f(+\infty)=1 \quad ({\rm CP}).
  \eeq
  Actually, we show that the solution of the CP has also  finite interface
  position $y_{0,{\rm CP}}>0$  that is larger than that for the above FBP.
  We then set $f(y) \equiv 1$ for $y \ge y_{0,{\rm CP}}$. With
  this continuation,
   the actual regularity of the CP solutions at
  the interface point $y=y_{0,{\rm CP}}$ is not known {\em a priori}, but the solutions
  are expected to be smoother than for the FBP,
  so (\ref{6}) are always valid. For instance, we show that,  for the TFE--4 with
  parameter $n>0$,
   \beq
 \label{n1}
  u_t= - (|1-u^2|^n u_{xxx})_x,
 \eeq
 the actual regularity of CP profiles at the interface
   is more
 than
 $C^{[\frac 3n]-1}$ ($[\cdot ]$ stands for the integer part), so
 solutions can be arbitrarily smooth for small $n>0$, \cite{Gl4}.

   Thus, as a key feature,
  in the CP,
   the profile $f(y)$ admits the trivial extension beyond the interface:
   $$
   f(y) \equiv 1 \quad \mbox{for} \quad y \ge y_{0,{\rm CP}}.
   $$
 For the FBP, such an extension is impossible, and actually, makes no sense, since,
 by definition, the FBP is  posed in a bounded moving  domain.
    Therefore, the problems with conditions (\ref{6}) and (\ref{7}) lead
to completely different similarity  profiles $f(y)$.

In general, we expect that, in the FBP setting, the RP--1
(\ref{4})--(\ref{6}) has a {\em countable} set of solutions
$\{f_k(y)\}$,
 where $f_k(y)$ approaches  the unique profile of the CP satisfying
 (\ref{7}), as $k \to \infty$. The first FBP profile is such that
 \beq
 \label{range1}
 f_1(y) \in (0,1) \quad \mbox{for all \,\,$y \in (0,y_0)$},
  \eeq
 and thus gives the stable (generic) solution of the Riemann problem for the TFE--4.
 Other FBP profiles $f_k(y)$ for all $k \ge 2$ have larger range, $\not \subset [0,1]$.
Though only the first $f_1(y)$ profile seems to make much physical
sense, we will study others to demonstrate mostly mathematical
properties of the ODE (\ref{4}).

A similar ``FBP--CP interaction" is shown to occur for the
RP--$\frac 12$.

\subsection{Connection with the linear bi-harmonic equation for $n=0$:
 towards oscillatory behaviour and branching in the CP}

It is curious that  structurally, the similarity solution of the
CP is close to the corresponding  solution for the linear
{bi-harmonic} equation
 \beq
 \label{l1}
 u_t=-u_{xxxx} \quad \mbox{in} \quad \re \times \re_+.
  \eeq
  In particular,
the {\em fundamental solution} of (\ref{l1}) has the similarity
form
 \beq
 \label{b1}
  \mbox{$
 b(x,t) = t^{-\frac 14} F(y), \quad y =  \frac x{t^{1/4}},
  $}
  \eeq
  where $F$ is the unique symmetric solution of the problem
   \beq
   \label{b2}
   \textstyle{
 {\bf B} F \equiv  - F^{(4)}+ \frac 14 \, F' y + \frac 14\, F=0 \quad \mbox{in}
   \quad \re, \quad \int F=1.
   }
 \eeq
The  kernel $F=F(|y|)$ is radial, has exponential decay,
oscillates as $|y| \to \infty$,
 and,   for some
positive constant $D$ and
  $d = 3 \cdot 2^{-11/3}$
 (see \cite[p.~46]{EidSys} and more
  details in Section \ref{SectBr}),
\beq
 \label{fbar}
 |F(y)| \le  D \,{\mathrm e}^{-d|y|^{4/3}}
\quad \mbox{in} \,\,\, \re.
 \eeq
Therefore,  by convoluting (\ref{b1}) with data (\ref{2}), the
following solution to the RP--1
 is obtained:
 \beq
 \label{b3}
  \mbox{$
u(x,t)= f_0(y) \equiv  b(x,t) * {\rm sign}\,  x=\Big(
\int\limits_{-\infty}^{y} -\int\limits^{+\infty}_{y} \Big) F(z) \,
{\mathrm d} z, \quad y=  \frac x{t^{1/4}}.
 $}
 \eeq
Hence, as earlier, $f(0)=f''(0)=0$, and $f(y)$ is oscillatory
about $\pm 1$ as $y \to \pm \infty$.

The linear PDE (\ref{l1}) is obtained from the more general TFE--4
 (\ref{n1})
 by passing  
  to the limit $n \to 0^+$. Such a continuous
 (homotopy) connection for some classes of solutions yields that
 the behaviour
 of solutions near $u = \pm 1$, which  for $n=0$ is oscillatory (there
 exist infinitely many intersections with $\pm 1$),   is
 inherited by solutions of the TFE--4 (\ref{n1}) at least for
 sufficiently small $n>0$. Actually, the oscillatory changing sign
 behaviour for TFEs takes place in the parameter interval \cite[\S~7.2]{Gl4}
  \beq
  \label{n12}
  n \in [0,n_h), \quad \mbox{with} \quad n_h=1.758665... \, ,
   \eeq
   where at $n=n_h$ there occurs a {\em heteroclinic bifurcation} in some
   related nonlinear ODE.


The existence and uniqueness of the solution $f_0(y)$ of Riemann's
problem for $n=0$ given in (\ref{b3}) suggests another approach to
the nonlinear one (\ref{4}). Namely, in Section \ref{SectBr}, we
intend to show that there exists a continuous $n$-branch of such
similarity profiles for the TFE--4 (\ref{n1}), which is originated
at $n=0$ from $f_0$.


\subsection{On extensions}

We must admit that Riemann's FBP and CP  for higher-order ODEs
(\ref{4}) are quite difficult, and we  can give a  rigorous
justification of existence for the FBP only.
 By a combination of
 analytic techniques
 and numerical evidence,
 we also treated the CP and establish that the similarity approach correctly describes the
 collapse of such  singularities for more general models.
  Furthermore, we can extend our
 similarity analysis
to $2m$th-order TFEs such as (the TFE--$2m$)
 \beq
 \label{8}
u_t =(-1)^{m+1} (|1-u^2|^n D_x^{2m-1} u)_x \quad (n > 0),
 \eeq
 with the discontinuous shock (\ref{2}),
for which we again pose
   the corresponding zero contact angle, zero-curvature, ...,  and zero-flux  FBP in a bounded domain and for
 the Cauchy problem in $\re \times \re_+$.
In both cases, the similarity solution for  the RP--1 is
 \beq
 \label{3m}
 u_+(x,t)= f(y), \quad \mbox{where} \quad y=\frac x{t^{1/{2m}}},
  \eeq
  and $f$ solves a  boundary-value problem for the ODE
  which, by anti-symmetry, we pose for $y>0$ only,
  \beq
  \label{4m}
    \textstyle{
   (-1)^{m+1}(|1-f^2|^n f^{(2m-1)})' + \frac 1{2m} \, f' y=0,
   \quad f(0)=f''(0)=... = f^{2m-2}(0)=0.
    }
     \eeq
 The zero contact angle FBP now involves $m+1$ ``zero contact angle" and flux conditions
 \beq
 \label{5m}
 f(y_0)=1, \quad f'(y_0)=f''(y_0)=...=f^{(m-1)}(y_0)=0, \quad
 (|1-f^2|^n f^{(2m-1)})(y_0)=0.
  \eeq

For the CP, we need just a single condition (\ref{7}).
 The case of the FBP is easier but still is a hard problem.
In the case of the CP for large $m$'s, any geometric-like shooting
techniques for the ODE (\ref{4m}) are illusive, since we have to
deal with truly  multi-parametric spaces. On the other hand, the
$n$-branching approach at the branching point $n=0$ (where there
exists the unique Riemann's profile $f_0(y)$ as in (\ref{b3}))
fits for any $m$ and implies existence of a suitable similarity
profile, at least, for $n>0$ small enough.

\smallskip

Concerning extensions of the fourth-order TFE (\ref{1}), we
consider the Cahn--Hilliard PDE with the extra unstable linear
diffusion term
 \beq
 \label{5m1}
u_t =- (|1-u^2| u_{xxx})_x - u_{xx}.
 \eeq
  However, this PDE is less physically motivated, though
  presents interesting interface properties.
Though the collapse for small $t>0$ of the initial singularity is
described by an analogous similarity solution (\ref{3}), we show
that the dynamical laws of interface propagation  changes
essentially in both the FBP and the CP. We derive the dynamic
interface equation in both cases, which establishes the connection
between the interface speed $s(t)$ and some interface differential
operators.

 Finally, we study the case of nonlinear unstable diffusion
 \beq
 \label{5m1N}
u_t =- (|1-u^2| u_{xxx})_x -
 (|1-u^2|^m u_x)_{x} \quad (m >0),
 \eeq
 which
 emphasizes some interesting
 mathematical phenomena.



\section{The FBP: on
 existence, uniqueness, analyticity, and stability}
\label{Sect2}

Here we consider the RP--1 (\ref{4})--(\ref{6}).  First, we notice
that the FBP and the CP have different profiles, and the latter
one has faster propagation with larger rescaled interface point,
 $$
 y_{0,{\rm FBP}} < y_{0,{\rm CP}}.
  $$



 \subsection{Local behaviour for the FBP: existence of a unique local smooth
 solution}

 We need to establish some basic local properties of the solutions
 near the free-boundary points. As we have mentioned, for the FBP
 we consider solutions in the proper range $f \in (-1,1)$, so that
 (\ref{4}) takes the form
  \beq
  \label{4N}
   \mbox{$
 -((1-f^2)f''')'+ \frac 14\, y f'=0.
  $}
  \eeq
  First of
 all, by classic ODE theory \cite[Ch.~I]{CodL}, $f(y)$ is smooth and analytic
at any non-degeneracy point, at which $f \not = \pm 1$. Close to
the interface at $y=y_0$, the local behaviour for solutions of
(\ref{4N}) is well-known \cite{BPW, FB0} and
 is given by the following expansion (still formal and hence to be
 justified):
 \beq
 \label{k1}
 \mbox{$
 f(y) = 1- C(y_0-y)^2 + \frac {y_0}{48} \, (y_0-y)^3 +
   \frac{ y_0 C}{960} \,  (y_0-y)^5+ O((y_0-y)^6),
  $}
 \eeq
 where $C > 0$ is an arbitrary fixed  constant.
 We will show that, actually, (\ref{k1}) is the expansion of an
 analytic function.
 Functions (\ref{k1})
 compose a 2D bundle with parameters $\{y_0,C\}$ to be matched
 with another 2D bundle at the origin
 with conditions (\ref{5}),
 \beq
 \label{k2}
 f(y)= A_1 y + A_3 y^3 + O(y^5),
  \eeq
  with constants $A_1>0$ and $A_3 \in \re$.
Actually, the asymptotics (\ref{k1}) give a smooth continuous
(relative to the presented parameters) penetration of orbits into
the
 space of analytic functions
  $f(y) \in (-1,1)$. Moreover, we will show that even (\ref{k1}) is
 the expansion of an analytic function at $y=y_0$, so we actually
 deal with globally analytic orbits.

 Since the papers \cite{BPW, FB0} dealt with
  third-order ODEs only obtained on integration once the divergent fourth-order counterpart
  (as a hint, see (\ref{tt2n}) below), we need to
  re-derive the expansion (\ref{k1}) for our fourth-order equation
  (\ref{4N}), to say nothing about analyticity, which, to our knowledge, was not
  addressed in the literature.
 It is curious that a standard application of contraction
 principle for the FBP then demand analysis in  functional spaces with
 singular (unbounded) weights.
  This  approach is rather general and applies to $2m$th-order ODEs for TFEs as in
  (\ref{4m}).

 \begin{proposition}
  \label{Pr.Exp}
 For any $C>0$, equation $(\ref{4N})$ has a unique local solution
 for $y<y_0$
 satisfying $(\ref{k1})$, which is strictly monotone increasing on
 some interval $y \in (y_0-\d,y_0)$.
  \end{proposition}

\noi{\em Proof.}
 We begin by integrating
 (\ref{4N}) over $(y,y_0)$, and using the free-boundary
 flux condition in (\ref{6})  obtain
  \beq
  \label{85}
  \mbox{$
  -((1-f^2)f''')(y) = \frac 14\, \int\limits_y^{y_0} f' y \,{\mathrm
  d}y.
   $}
   \eeq

In order to simplify further calculus,
 let us explain the origin of the unique monotone solution
 of (\ref{85}). We represent the right-hand side in the equivalent
 form
 \beq
  \label{851}
  \begin{matrix}
  -((1-f^2)f''')(y) = \frac 14\, \int_y^{y_0} f' y \,{\mathrm d}y
  \equiv - \frac 14 \,  \int_y^{y_0}[(1-f)'y]\,{\mathrm d}y\qquad
  \smallskip\smallskip\smallskip\\
  = - \frac 14\, \big\{ \int_y^{y_0}[(1-f)y]'\,{\mathrm d}y - \int_y^{y_0}(1-f)\,{\mathrm
  d}y\big\}
   = \frac 14\,(1-f)y +  \frac 14\, \int_y^{y_0}(1-f)\,{\mathrm
   d}y. \qquad
   \end{matrix}
   \eeq
Therefore, dividing by $1-f^2=(1-f)(1+f)$ yields the following
perturbed equation:
 \beq
 \label{882nn}
 \mbox{$
 - f'''= h(f,y)=h_1+h_2 \equiv \frac y{4(1+f)}+ \frac 1{4(1-f)(1+f)}\,  \int_y^{y_0}(1-f)\,{\mathrm
 d}y.
  $}
  \eeq
 Obviously, the first term is non-singular at the point $(y_0,1)$ for $f \approx 1$,
 $y \approx y_0$,
 so, by classic theory,
 the unperturbed ODE
 \beq
 \label{882nns}
 \mbox{$
 - f'''= \frac y{4(1+f)}
  $}
  \eeq
 possesses a
  unique local {\em analytic} solution
 satisfying,
 according to the asymptotics (\ref{k1}):
  \beq
  \label{883nn}
   \mbox{$
  f(y_0)=1, \quad f'(y_0)=0, \quad f''(y_0)=-2C<0, \quad
  f'''(y_0)=- \frac {y_0}8<0.
  $}
  \eeq
 Note that any $C^3$-solution satisfying (\ref{883nn}) is strictly monotone
 increasing for $y \approx y_0^-$.

 The second term in (\ref{882nn})
 contains certain weak singularity
 in the same class of smooth solutions, so we should treat it in
 a weighted space.
 As a standard procedure of Perron--Picard--Lyapunov type, we next integrate
 (\ref{882nn}) three times with conditions (\ref{883nn}) to obtain
 the equivalent integral equation
  \beq
  \label{int1}
 \mbox{$
  f(y)= 1- C(y_0-y)^2 - \int \int
  \int_y^{y_0} h(f(y),y)\,({\mathrm d} y)^3.
   $}
   \eeq
It is convenient to rewrite (\ref{int1}) for the functions $v(y)$
given by shifting
 \beq
 \label{fff1}
  \mbox{$
 f(y)= 1- C(y_0-y)^2 + v(y).
 $}
 \eeq
 Then the integral equation reads
  \beq
  \label{fff1I}
 \mbox{$
  v(y)=  {\mathcal M}(v)(y) \equiv - \int \int
  \int_y^{y_0} \hat h(v(y),y)\,({\mathrm d} y)^3,
   $}
   \eeq
 where we denote $\hat h(v,y)= h(1- C(y_0-y)^2 + v,y)$.

We next introduce a functional framework,
 which is suitable for
 integral equation with  singularities; see
 \cite[\S~5.8]{AMGV}, where these techniques were applied to
 uniform stability problems for degenerate Hamilton--Jacobi and
singular perturbed parabolic PDEs. Thus, let us
  fix a small $\d>0$ and set $I_\d=[y_0-\d,y_0)$. We consider the
integral equation  (\ref{fff1I}) in the functional space
 \beq
 \label{x1}
  \mbox{$
 X_\rho= \{v \in C(I_\d), \,\, v(y_0)=v'(y_0)=0, \,\, \rho v \in
 L^\iy(I_\d)
 \}, \,\,\, \mbox{where} \,\,\,  \rho=  \frac
 1{(y_0-y)^3}.
 $}
  \eeq
We endow  $X_\rho$  with the following  sup-norm:
 \beq
 \label{ss1n}
  \mbox{$
  | v |_\rho= \sup_{y \in I_\d} \, \{\rho(y) \,
  |v(y)|\},
   $}
   \eeq
   i.e., we  introduce the singular weight $\frac
   1{(y_0-y)^3} \to \iy$ as $y \to y_0^-$ into the standard
   framework of the space $C(I_\d)$ of continuous
   functions with necessary boundary conditions at $y_0=0$.
    With this distance
(\ref{ss1n}), $X_\rho$ becomes a complete metric space being a
closed subspace of the Banach space $C_\rho(I_\d)$ of twice
differentiable functions at $y=y_0$,
 for which the semi-norm
 \beq
 \label{x2}
 \mbox{$
 |v|_\rho =\sup_{y \in I_\d} \, \big\{\rho(y) \,
  \big|v(y)-v(y_0)-v'(y_0)(y-y_0) - \frac 12\, v''(y_0)(y- y_0)^2\big|\big\}
   $}
   \eeq
 is finite. The space $C_\rho(I_\d)$ has the natural norm
  $$
  \|v\|_\rho= |v(y_0)|+|v'(y_0)|+ |v''(y_0)|+ |v|_\rho.
   $$

Firstly, it is an easy exercise to see that two terms in
(\ref{882nn}) being rewritten for $v$, immediately yield
 \beq
 \label{int4}
  {\mathcal M}:X_\rho \to X_\rho \quad \mbox{for all small
  $\d>0$}.
  \eeq
Indeed, (\ref{882nn}) implies that, for $ y \approx y_0$,
 $
\hat h \sim \frac {y_0}8 + O((y_0-y)),
  $
 so, on triple integration in (\ref{fff1I}), we obtain the result:
  $$
  \mbox{$
 0 \le  {\mathcal M}(v)  \le \big[\frac {y_0}8+ o(1)\big]
 (y_0-y)^3 \le y_0 (y_0-y)^3 \,\,\, \mbox{in} \,\,\, I_\d.
 $}
 $$

Secondly, checking the contractivity of ${\mathcal M}$ on
$X_\rho$, one can see that the first term $h_1$ in (\ref{882nn})
does this as inducing non-singular and analytic operator. Consider
the second term $h_2$ therein.
  Taking
arbitrary $v_{1,2} \in X_ \rho$ and denoting by $f_{1,2}$ their
$f$-counterparts due to the change in (\ref{fff1}), setting for
convenience $\D v= v_1-v_2 = f_1-f_2$ will provide us
 by standard manipulations similar to the
Lagrange formula of finite increments  the following:
 \beq
 \label{int5}
  \begin{matrix}
 h_2(f_1)- h_2(f_2)= \hat h_2(v_1)-\hat h_2(v_2) \ssk\ssk
 \\
 =
  \frac y{4(1+f_1)} + \frac { \int(1-f_1) }{4(1-f_1^2)} -
\frac y{4(1+f_2)} - \frac {\int(1-f_2)}{4(1-f_2^2)} \ssk\ssk
 \\
 =  \frac {y (f_2-f_1)}{(1+f_1)(1+f_2)}
+ \frac { \int(1-f_1)}{4(1-f_1^2)} - \frac
{\int(1-f_2)}{4(1-f_1^2)}
 + \frac {\int(1-f_2)}{4(1-f_1^2)}  - \frac
{ \int(1-f_2)}{4(1-f_2^2)}\ssk\ssk \\
   = -\frac {y
\D v}{(1+f_1)(1+f_2)}- \frac {\int \D v}{4(1-f_1^2)} +
 \frac { \int
 (1-f_2)}4 \big( \frac 1{1-f_1^2} - \frac 1{1-f_2^2} \big).
 \end{matrix}
 \eeq
 Obviously, for small $\d>0$, we may assume that, in $I_\d$,
  \beq
  \label{as1}
   \begin{matrix}
    1 \le 1+f_{1,2} \le 2, \,\,\,
     \frac 12 \le \frac 1{1+f_{1,2}} \le 1,\qquad\qquad
 \ssk\ssk\\
 \frac C2(y_0-y)^2 \le 1-f_{1,2} \le 2C(y_0-y)^2, \,\,\,
 \frac 2{C(y_0-y)^2}
 \le  \frac 1{(1-f_{1,2}^2)} \le \frac 4{C(y_0-y)^2}.\qquad\qquad
   \end{matrix}
   \eeq
 Then, denoting
 here and later on by  $A>0$ various constants that are independent of
 $\d$,
 we obtain the following estimate:
 \beq
 \label{int50}
 \begin{matrix}
|\hat h_2(v_1)- \hat h_2(v_2)| \le \frac{y_0}4 | \D v| +
 \frac 1{C(y-y_0)^2} \, \int|\D v| + \frac 14\, \frac{2C(y_0-y)^3}3 \,\,
 \frac {2|\D v|}{
 \frac{C^2}4\, (y_0-y)^4} \qquad\qquad
  \\
  \le
 A \big[\, |\D v|  + \frac 1{(y_0-y)^2} \, \int |\D v| +
  \frac{|\D v|}{y_0-y}\, \big].\qquad\qquad
   \end{matrix}
 \eeq
Using the fact that, in the metric (\ref{ss1n}),
 \beq
  \label{dddd1}
 |\D v(y)| \le |\D v|_\rho (y_0-y)^3,
 \eeq
 we obtain that
  \beq
  \label{int51}
|\hat h_2(v_1)- \hat h_2(v_2)| \le A | \D v|_\rho \big[(y_0-y)^3 +
(y_0-y)^2 \big] \le  A (y_0-y)^2 | \D v|_\rho.
 \eeq

Finally, using the metric (\ref{ss1n}) again and substituting
(\ref{int51}) into (\ref{fff1I}), we then obtain the
 following principal bound:
 \beq
 \label{ss4}
 \begin{matrix}
  | {\mathcal M}(v_1)- {\mathcal M}(v_2) |_\rho
   \le
     \rho(y) \, \int \int \int_{y}^{y_0}  A(y_0-y)^2 \, | \D
     v|_\rho \, ({\mathrm d}y)^3
     \smallskip
\smallskip\smallskip\\
     \le | \D v |_\rho A  \frac 1{(y_0-y)^3} \,
\int \int \int_{y}^{y_0}  (y_0-y)^2 \, ({\mathrm d}y)^3
 \le
\frac A{60} \, \d^2 \,\,  | \D v |_\rho .
  \end{matrix}
  \eeq
This implies that
 \beq
 \label{int6}
  \mbox{${\mathcal M}$ \,\,\, is a contraction in
  $X_\rho$ for all small $\d>0$ such that $\frac A{60} \, \d^2<1$ in (\ref{ss4})}.
  \eeq
   Hence,
 by Banach's Contraction Principle
(see e.g. \cite[p.~206]{KrasZ}), for all sufficiently small
$\d>0$,  (\ref{int4}) and (\ref{int6}) guarantee
 existence and uniqueness of a suitable local
  solution of (\ref{fff1I}) in $X_\rho$ for any fixed $C>0$.
 $\qed$

\smallskip



\subsection{Analyticity at $y=y_0$}

This is also a principal question of TFE theory. We begin with two
simpler illustrations.

\ssk

\noi\underline{\em 1. Explicit analyticity for the standard TFE}.
  Consider
    the ``standard" TFE
 \beq
 \label{tt1nn}
 u_t=-(u u_{xxx})_x \quad \mbox{in}
 \quad \re \times \re_+,
  \eeq
 which has the same type of degeneracy as (\ref{eq1}) but now at
the equilibrium $\{u=0\}$, so that equations near degeneracy sets
are approximately connected by the change $u \mapsto 1-u$.
  Looking for the source-type solution of (\ref{tt1nn}) leads to a simpler
  ODE in divergence form:
   \beq
   \label{tt2}
   \mbox{$
   u_{\rm s}(x,t)=t^{-\frac 1{5}} f(y), \quad y=\frac x{t^{1/5}}
   \,\, \Longrightarrow \,\,
   -(f f''')'+ \frac 1{5}\,(f y)'=0.
    $}
    \eeq
The ODE in (\ref{tt2}) near the interface is again approximately
connected with (\ref{4N}) by $f \mapsto 1-f$.
 Integrating (\ref{tt2}) ones and dividing by $f>0$ yields the
following explicit  solution detected  by Smyth and Hill  in 1988
\cite{Smyth88}:
 \beq
   \label{tt2n}
   \mbox{$
   - f f'''  +\frac 1{5}\,f y=0 \,\, \Longrightarrow
   \,\, f'''=  \frac 1{5}\, y, \,\,\,\, \mbox{i.e.,}\,\,\,\,
   f(y) =  \frac {1}{{120}}\,(y_0^2- y^2)^2.
    $}
   \eeq
Thus, we see a finite polynomial solution (\ref{tt2n}), which is a
perfect and simplest analytic solution at the interface $y=y_0$,
which, surely, is a singularity degeneracy point for the ODE in
(\ref{tt2}). As an easy  illustration, we observe that the
interface at $y=y_0$ occurs via the geometric configuration when
 \beq
 \label{dd1ss}
 \mbox{an analytic function $f(y)$ ``touches down" the level
 $\{f=0\}$ at some $y=y_0>0$}.
  \eeq
 Being formally
extended by $f(y) \equiv
 0$ for $|y| \ge y_0$ (in fact, this makes no particular sense for the FBP posed in $\{|y| < y_0\}$
 only),
  the solution (\ref{tt2})  solves another singular
``Riemann
 problem" with a measure as initial data,
 $$
  \mbox{$
 u_{\rm s}(x,t) \to  c_0 \, d(x)
 \quad \mbox{as} \quad t \to 0^+ \quad \big(\, c_0= \frac {2 y_0^5}{15}\,
 \big)
  $}
 $$
 in the sense of distributions.
 Hence, (\ref{tt2}) is
 a formal {\em fundamental solution} of the FBP.

Unfortunately, as we have demonstrated, our ODE (\ref{4N}) is more
complicated and does not admit an explicit integration. However,
we are going to show that, for the FBP,
 analyticity of the similarity profile $f(y)$ as in (\ref{tt2n}) is
 not exceptional and exclusive, so
 $f(y)$ of (\ref{4N}) remains
analytic and actually is close to that in (\ref{tt2n})  as $y \to
y_0$, i.e., with the change $f \mapsto 1-f$, the FBP profiles for
(\ref{4N}) for $y \approx y_0$
 \beq
 \label{ana1ss}
 \mbox{
 are given by small analytic perturbations of the analytic
 solution in (\ref{tt2n}).}
  \eeq
Therefore, the explicit finite polynomial FBP solution in
(\ref{tt2n}) is a good illustration of the origin of the real
analyticity of other non-explicitly given solutions. In its turn,
the analyticity has also a strongly motivated natural origin in
complex variable presentation:

\ssk

\noi\underline{\em 2. Analyticity for complex-valued ODEs and
Weierstrass theorem}. The origin of the analyticity of the FBP
profiles can be also seen from the ODE theory in the complex plane
$z=x+{\rm i} y \in {\mathbb C}$, where $x$ now stands for $y$ in
(\ref{4N}). Namely, consider an analytic function
 \beq
 \label{ff1ss}
 f=f(z)=u(x,y)+{\rm i} v(x,y) \quad (u_x=v_y, \,\,\, u_y=-v_x)
  \eeq
  satisfying the complex ODE
   \beq
   \label{ff2}
    \mbox{$
   -((1-f^2)f_{zzz})_z + \frac 14 \, z f_z=0 \quad \mbox{in}
   \quad {\mathbb C}.
    $}
    \eeq

Since (\ref{ff2}) contains analytic nonlinearities only, it admits
a family of analytic solutions in a neighbourhood $O$ of the
origin $z=0$
 \beq
 \label{zz1}
  \mbox{$
 f(z)= \sum_{(k \ge 0)} c_k z^k,
  $}
  \eeq
 with the following relation on the
expansion coefficients:
 \beq
 \label{zz1n}
   \begin{matrix}
 c_{k+4}(k+1)(k+2)(k+3)(k+4) \ssk\ssk\\ - \sum\limits_{l\ge 3,0
\le  j \le k+4-l}
 c_l c_{k+4-l-j}c_j \, l(l-1)(l-2)(k+1)= \frac 14 \, k c_k, \,\,\, k
 \ge 0.
  \end{matrix}
 \eeq
 Choosing an analytic solution that is better associated with the
 FBP expansion (\ref{k1}), i.e.,  $c_0=1$ and $c_1=0$,
 \beq
 \label{zz1N}
  \mbox{$
 f(z)=1+ \sum_{(k \ge 2)} c_k (z-z_0)^k, \quad \mbox{where}
 \quad z_0=(x_0,0),
  $}
  \eeq
 yields a similar algebraic system,
 \beq
 \label{zz1nN}
   \begin{matrix}
- \sum\limits_{l \ge3,\, 1 \le j \le k+4-l} c_j c_l c_{k+4-l-j}\,
 l(l-1)(l-2)(k+1) \ssk\ssk\\ = \,
 \frac 14 \, k c_k + \frac 14\, x_0(k+1)c_{k+1}, \quad k
 \ge 1.
  \end{matrix}
 \eeq
 In fact, the analyticity of the solutions, i.e., non-zero radius
 of convergence of the series in (\ref{zz1}) and ({\ref{zz1N}) can
 be seen from the infinite algebraic systems (\ref{zz1n}) and
 (\ref{zz1nN}). In other words, one needs to show that
 the recurrent relations generated by (\ref{zz1n}) or (\ref{zz1nN}) admit solutions
 $\{c_k\}$ without huge growth such as
   \beq
   \label{hh1}
    \mbox{$
    \limsup_{k \to \iy} |c_k|^{\frac 1k}=\iy \quad (\mbox{then}
    \,\,\, R=0).
     $}
    \eeq
 A direct proving that (\ref{hh1}) does not
 take place for the above power series is not straightforward (but indeed seems doable),
 so we
 will follow a simpler  approach.

Meantime, we note that,
 on the other hand, iterating the integral equation equivalent to (\ref{ff2}) via the simple
  iteration
starting with an analytic initial data $f_0(z)$ in $O$
(this will be done shortly for the equivalent integral equation in
the real case, so we do not discuss details)
 will give us a sequence of analytic functions
$\{f_k(z)\}$. In case of the uniform convergence of the series
along a subsequence in a neighbourhood of $z=0$ (again, the
convergence  is naturally proved for the integral counterpart of
the equation),
 $$
 f_k(z) \to f(z) \quad \mbox{as} \quad k=k_j \to \iy,
 $$
 will provide us with an analytic solution, as a manifestation of
 Weierstrass' classic theorem (the limit of a uniformly converging sequence of
 analytic functions in $O$ is also analytic).

 This is the origin of analytic solutions of the degenerate ODE
 (\ref{ff2}) in the complex plane. Let us trace out the link to the
 real case. In the variables in (\ref{ff1ss}), currently using
 the differentiation $f_z= u_x+{\rm i}v_x$,
  the ODE
 (\ref{ff2}) reduces to a system for
 ${\rm Re}\,f$ and ${\rm Im}\, f$,
  \beq
  \label{zz3}
   \left\{
    \begin{matrix}
    -[(1-u^2+v^2)u'''+ 2 uv v''']' + \frac 14 \, (x u'-y v')=0,
    \ssk\ssk\\
 -[(1-u^2+v^2)v'''- 2 uv u''']' + \frac 14 \, (y u'+x v')=0,
  \end{matrix}
   \right.
    \eeq
    where $'=D_x$. It then follows that:
  \beq
  \label{zz4}
  \mbox{(\ref{4N}) is (\ref{ff2}) at the axis $\{y=0\}$, on which
  $v(x,0)=0$.}
   \eeq
    Indeed, then also $v_x(x,0)=v_{xx}(x,0)=... \equiv 0$
    that make the second equation in (\ref{zz3}) tautological, and
    then the first one for $f(x)=u(x,0)$ coincides with  (\ref{4N}), where $y \mapsto x$.
 One can see that an analytic solution of (\ref{ff2}) satisfying
 (\ref{zz4}) must have the form (\ref{zz1}) with all {\em real}
 expansion
 coefficients\footnote{It is convenient to call such solutions
 R-analytic, with
 ``R" standing for ``Real coefficients" (not to confuse with ``real analytic" to be
 treated below).} $\{c_k\}$.
 Then, among those analytic solutions, there can be some, which we
 have introduced in
 (\ref{zz1N}),
 corresponding to a ``touching down configuration" of $u(x,0)$ as
 in (\ref{dd1ss}).
A proper solvability of the analytic elliptic system (\ref{zz3})
with the condition at the $x$-axis  in (\ref{zz4}), i.e., in the
class of R-analytic solutions, and with the geometry as in
(\ref{dd1ss}) is not
  a part of the present study. However, this
shows
 a key analyticity link between complex and real ODEs and the role
 of Weierstrass' theorem in constructing analytic solutions of such degenerate
 ODEs with analytic nonlinearities.

\ssk

\noi\underline{\em 3. Main result}. Thus, we now turn our
attention to the real valued ODE (\ref{4N}) and prove:

\begin{proposition}
 \label{Pr.An}
 The local solution
 of Proposition $\ref{Pr.Exp}$ is real analytic at $y=y_0$.
  \end{proposition}

\noi{\em Proof.} \underline{\em Uniqueness.}
 First, let us note again
 that in (\ref{882nn}), the main
  first non-singular term itself produces  existence of an analytic solution.
 At the same time, the second perturbation term is also analytic
in the sense that it maps admissible (i.e., those that do not
create singularities in the integral) analytic functions into
analytic. Therefore, the ODE (\ref{882nn}) admits a unique formal
analytic expansion.
 Indeed, setting
   \beq
   \label{aa1}
  \mbox{$
  f(y)= 1+ \sum_{(j \ge 2)} c_j \D^j, \quad \mbox{where} \quad
 \D=y-y_0 \quad (c_2=-C < 0)
 $}
 \eeq
 and substituting this formal expansion into (\ref{882nn})
  yields another expansion
 \beq
 \label{aa2}
  \mbox{$
  h(f(y),y) \equiv \frac 1{4(1- f)(1+ f)}
  \int_y^{y_0} f= \sum_{(k \ge 1)} b_k \D^k,
  $}
   \eeq
 with a uniquely solvable relation between
 coefficients $\{c_k\}$ and $\{b_k\}$ as in (\ref{zz1n}). Integrating (\ref{882nn}) three times
 leading to (\ref{fff1I}) preserves this formal
  power expansion.
 Thus,  this yields the uniqueness:
 if an analytic solution (\ref{aa1}) of (\ref{fff1I}) exists, it
 is unique.

 Note also that  this implies that the solution constructed in
 Proposition \ref{Pr.Exp} is at least $C^\iy$, i.e., applying
 finite expansion analysis with the remainders of order
 $O(\D^l)$, with arbitrarily large $l \gg 1$, yields that there exist all the derivatives
  \beq
  \label{all1}
    \exists \,\, f^{(j)}(y_0) \equiv j! \, c_j \quad \mbox{for
    all} \quad j \ge 0.
     \eeq
Of course, this is easily seen from (\ref{882nn}) by
differentiating as many times as necessary:
 $$
  \mbox{$
 f''(y_0)=-2C, \quad f'''(y_0)=
 \frac {y_0}8,
 \quad f^{(4)}(y_0)=0, \quad f^{(5)}(y_0)=- \frac {y_0C}{960},\,\,
 ...\, ,
  $}
  $$
 which lead to the expansion (\ref{k1}).

\ssk

\noi\underline{\em Existence}. Let us show
convergence of the series (\ref{aa1}). First of all, there exists
an {\em a priori} bound for (\ref{882nn}): for any suitable
solutions $f$ such that $v \in X_\rho$, in $I_\d$,
 \beq
 \label{aa3}
  \mbox{$
 0 < h(f,y) \le  \frac {y_0}4 +  \frac 1{2C(y_0-y)^2} \int 2C
 (y_0-y)^2 \,{\mathrm d}y \le \frac {y_0}4 + \frac {y_0-y}3 \le
 y_0.
  $}
   \eeq

We now perform an analytic simple iteration of the integral
equation (\ref{fff1I}),
 \beq
 \label{ana1}
 v_{k+1}= { \mathcal M}(v_k), \quad k=0,1,2,... \, ,
 \eeq
 by taking analytic initial data
  \beq
  \label{om1}
 v_0=0 \in X_\rho, \quad \mbox{i.e.,} \quad f_0(y)=1-C(y_0-y)^2.
  \eeq
 By Banach's Contraction Principle, we then obtain
 a sequence of functions $\{f_k(y)\}$ such that
 the corresponding  polynomial sequence $\{v_k(y)\}$
 given by (\ref{fff1}) satisfies as $k \to \iy$
 \beq
 \label{ana2}
 v_k(y) \to v(y) \quad \mbox{in $X_\rho$ and uniformly in
 $I_\d$}.
  \eeq
 These functions are given by converging power series
  \beq
  \label{tk1}
   \mbox{$
   f_k(y) =1+ \sum_{(j\ge 2)} c_j^{k} \D^j.
    $}
    \eeq
    By (\ref{aa3}), the  polynomial sequence
    $\{f_k(y)\}$ is uniformly bounded in $I_\d$.

Finally, we then need to prove that the radiuses of convergence
$R_k$ of the approximating series (\ref{tk1}) are bounded from
below:
 \beq
 \label{rk1}
 R_k \ge \d_0>0 \quad \mbox{for all} \quad k \gg 1.
  \eeq
  To this end we need some more calculus concerning the integral
  operator involved.

 Thus, according to
(\ref{int1}), by the change
 $$
\zeta=y_0-y>0, \quad  f(y)= 1-C \zeta^2(1+\e(\zeta))
 $$
  equation
(\ref{882nn}) reduces to
 \beq
 \label{yy1}
  \mbox{$
   C(\zeta^2 \e)'''=  \bar h (\e,\zeta)
\equiv  \frac {y_0-\zeta
  }{4(1+f)} + \frac 1{4(1+f) (1+\e)} \, \frac 1{\zeta^2}\, \int_0^\zeta \zeta^2
  (1+\e) {\mathrm d} \zeta,
   $}
   \eeq
where now the right-hand sided is  analytic
 in $\e$ but is not regular in $\zeta$ in the last term with
 formal singularity $\frac 1{\zeta^2}$, which we have to pay the
 main and the only our attention.
  For simplicity and without fear of confusion, we still keep $f$ in the terms that are obviously
regular analytic and provide no difficulties. Then changing the
independent variable in the left-hand side
 \beq
 \label{yy2}
  \mbox{$
  \zeta=\phi(s)>0 \,\,\Longrightarrow \,\,
  (\zeta^2 \e)'''= \zeta^2 \e'''+...=
  \frac{\phi^2}{(\phi')^3} \, \e'''+...
 $}
  \eeq
 we set $\phi'= \phi^{2/3}$, i.e., $\zeta=\phi(s)= (\frac s3)^3$ gives an analytic change of
  variables,
  to
 obtain
  \beq
  \label{yy3}
   \mbox{$
  C \e'''+...= \bar h\big(\e,(\frac s3)^3\big)=...+ \frac
  1{4(1+f)(1+\e)} \,  \frac {1}{9 s^6} \int_0^s s^8(1+\e)\, {\mathrm
  d}s,
   $}
   \eeq
   where we omit lower-order regular terms in $s$, which supply
   us with guaranteed analytic $s$-expansions.
   Actually, the integral representation (\ref{yy3}) actually completes
 the proof of analyticity. Indeed, we already know that the
 solutions is unique, and then substitution into (\ref{yy3}) the
 analytic series gives a rather standard relation between
 expansion coefficients that guarantee the non-zero radius of
 convergence. Let us present some comments.

Consider the iteration (\ref{ana1}) in the  $\e$-variable that
gives us the required sequence $\{\e_k,\, k \ge 0\}$,
 which thus is given by the recursion
\beq
  \label{yy3N}
   \mbox{$
  C \e_{k+1}'''+...= \bar h\big(\e_k,(\frac s3)^3\big)=...+ \frac
  1{4(1+f_k)(1+\e_k)} \,  \frac {1}{9 s^6} \int_0^s s^8(1+\e_k)\, {\mathrm
  d}s.
   $}
   \eeq
 Let us perform the first step to estimate the radius $R_1>0$.
 Let
$C=y_0=1$ and take
 $$
  \mbox{$
 \e_0=0 \quad \mbox{and} \quad f_0(\zeta)=1-\zeta^2= 1- \big(\frac s3
 \big)^6.
  $}
 $$
This yields that $\e_1(s)$ is given by the triple integration of
the equation
 $$
  \mbox{$
 \e_1 '''+...=...+ \frac 1{2-(\frac s3)^6} \frac 1{s^6}
 \int_0^s s^8 \, {\mathrm d}s,
 $}
  $$
  where we again keep the ``non-analytic" term only that, by assumption, may
  eventually vanish the radius of convergence of approximations $\e_k(s)$.
We also omit all the constants that do not affect this radius.
Thus, $\e_1$ is given by integrating
 \beq
 \label{pp1}
  \mbox{$
  \e_1'''+...=...+ \frac {s^3}{2-(\frac s3)^6}
  \quad \Longrightarrow \quad R_1^{\rm n-a}=3 \cdot 2^{\frac 16} \ge 1,
   $}
    \eeq
  where we indicate the radius of convergence for the power series
  for $\e_1(s)$, which is determined by the formally non-analytic term
  (the actual radius can be smaller if the analytic terms define a
  smaller radios, that does not happen). Thus, $\e_1(s)$ admits the
  unique analytic continuation denoted again by $\e_1(z)$ into the disc
  $B_1=\{|z|<R_1\} \subset {\mathbb C}$.


Consider the transition $k \mapsto k+1$ via (\ref{yy3N}).
 Let $\e_k(z)$ be analytic in $B_k=\{|z| < R_k\}$, where $R_k >0$.
By the definition of $\e$, we may assume by (\ref{ana2}) that
 \beq
 \label{ee1J}
  \mbox{$
  \e_k(0)=0 \quad \mbox{and}
  \quad
 |\e_k(s)| \le  \frac 12
 \quad \mbox{in}
 \quad I_\d \quad(\mbox{actually} \,\,\,  |\e_k(s)| \le |s|),
  $}
  \eeq
  and moreover we also may assume that the analytic continuation
  $\e_k(z)$ also satisfies this estimate in $B_k$.
 Consider the right-hand side of (\ref{yy3N}). Obviously both
 functions $1+f_k$ and $1+\e_k$
 are analytic in $B_k$, and so do $\frac 1{1+f_k}$ and $\frac
 1{1+\e_k}$ provided that these do not vanish in $B_k$ that is
 true by (\ref{ee1J}) for small $\d>0$. The integration in
 (\ref{yy3N}) also does not change the radius of convergence.
 Then, overall, the product of three analytic functions in
 (\ref{yy3N}) preserves the same radius of convergence $R_k$, so
 that by triple integration $\e_{k+1}(s)$ will inherit at least the same
 uniform estimate of $R_{k+1}^{\rm n-a}$ from below.


Thus, inside this interval of uniform convergence of all the
series, we get the desired result:
  as $k \to \iy$, uniformly
in $I_\d$,
 \beq
 \label{tk4}
  \mbox{$
  f_k(y)=1+ \sum_{(j \ge 2)} c_j^k \D^j
 \to 1+ \sum_{(j \ge 2)} c_j \D^j = f(y),
    $}
    \eeq
  where $f(y)$ is the unique solution from Proposition \ref{Pr.Exp}.
  This implies that $f(y)$ is analytic at the degeneracy point $y=y_0$.
     $\qed$

\ssk

It is worth mentioning that in the above analyticity proof, we
 essentially used the existence-uniqueness via contractions in
 $C_\rho$
 established in Proposition \ref{Pr.Exp}.
 On the other hand, the presented estimates of the radiuses of
 convergence of $\{\e_k(z)\}$ can also provide existence and
 uniqueness of analytic solutions, but this issue should be
 pre-converted into the complex variable nature for the ODE
 (\ref{ff2}). Then, this analysis will precisely emphasize
the application of the already announced Weierstrass theorem.
However, proving non-existence of  non-analytic solution profiles
will anyway require a contractivity study of the integral equation
such as (\ref{fff1I}) in metric spaces of functions of finite
regularity.

\subsection{First similarity profile: numerics}

Concerning  existence of finite interface in the degenerate ODE
such as (\ref{4}), we refer to first results in \cite{Bern88,
BMc91}.
Existence (and uniqueness) of the first similarity profile
$f=f_1(y) \in [0,1]$ depends on a  2D--2D matching argument,
where, using two parameters in the bundle (\ref{k2}), we need to
satisfy three conditions (\ref{6}) with the unknown interface
point $y_0$ (a free parameter); see below.



Shown in Figure \ref{F01} are results of numerical shooting for
the FBP (\ref{4})--(\ref{6}). We indicate here a number of
profiles (dotted lines) satisfying the anti-symmetry conditions at
the origin (\ref{5}) and two conditions at the right-end point
 \beq
 \label{76}
 f(y_0)=1 \quad \mbox{and} \quad f'''(y_0)=0.
  \eeq
We next vary the length $y_0>0$ to get the zero contact angle
 \beq
 \label{77}
 f'(y_0)=0,
 \eeq
 and then, by construction, all three conditions (\ref{6}) are
 valid.
We  obtain the first profile $f_1(y)$ (the boldface line)
satisfying the necessary range condition (\ref{range1}), which was
numerically clearly unique, with the interface at
  \beq
   \label{ii1}
 y_{01} \sim 2.35.
  \eeq
  The above shooting procedure also  indicates that there
  exists the unique second FBP profile $f_2(y)$ with the interface at
 $$
 y_{02} \sim 3.75.
 $$
We return to  this   multiplicity problem of the FBP profiles in
Section \ref{Sect31} after introducing and studying the unique
oscillatory similarity profiles of the Cauchy problem, to which
FBP ones turn out to have a direct relation.



\begin{figure}
\centering
\includegraphics[scale=0.85]{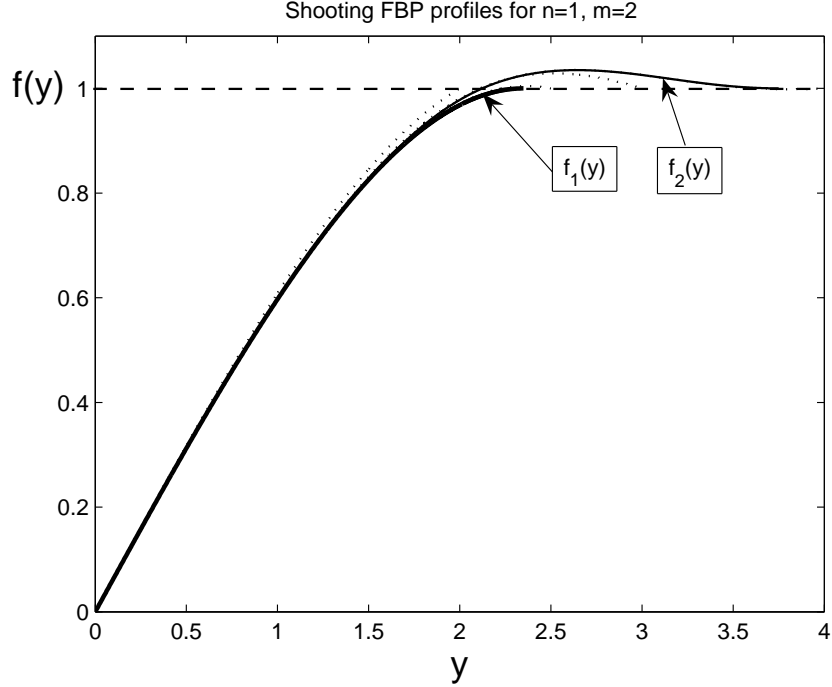} 
\vskip -.3cm \caption{\small Construction of the unique first
similarity profile $f_1(y)$ for the FBP.}
   \vskip -.3cm
 \label{F01}
\end{figure}

 \subsection{On existence and uniqueness of the first monotone FBP $f_1$}

We begin by mentioning that the mathematical techniques
 developed in Bernis \cite{Bern88} and Bernis--McLeod \cite{BMc91}
  as well as in \cite{BPW, FB0},
  are essentially directed to
third-order ODEs obtained on integration ones via conservation
laws, so these do not directly apply to the truly fourth-order ODE
(\ref{4}). Nevertheless,  the methods and results therein can
 be
 useful.
In particular, Bernis energy estimates \cite[\S~7]{BMc91} settles
finite propagation property for such ODEs.

 Shown in Figure \ref{F01} is the actual mathematical strategy to
 obtain
 the necessary profile $f_1(y)$ to get the following:

\begin{proposition}
 \label{Pr.11}
  The problem $(\ref{4})$--$(\ref{6})$ admits a solution
  $f_1(y)$ satisfying
   \beq
   \label{95}
   f_1'(y) >0, \,\,\, f_1''<0 \quad \mbox{on}
   \quad (0,y_0), \quad \mbox{so}
   \quad 0 < f_1(y)<1.
 \eeq
 \end{proposition}

\noi{\em Proof.} Thus,
  we shoot from the
 unknown interface point $y=y_0>0$ according to the bundle
 (\ref{k1}), and denote the solution as
  \beq
  \label{91}
  f=f(y;y_0,C),
  \eeq
  which can be continued until the regular point $y=0$.
  One can see that the ODE (\ref{4}) does not admit local
  singularities at finite points.
 Our further analysis uses the  continuous dependence of $f(y)$ on the parameters $y_0$ and
 $C$. Note that, as we have shown,
  the dependence is also continuous in the expansion
  (\ref{k1}) near singular (analytic) points.

  We next need a few  simple observations concerning such shooting
  orbits:

 (i) We want to know that an inflection point of $f(y)$, where
 \beq
 \label{94}
 f''(y_1)=0, \quad y_1 \in [0,y_0), \quad f(y_1) \in [0,1),
 \eeq
 actually occurs.
  This is seen from the ODE and, moreover, from
 the expansion (\ref{k1}) that is obtained from it.
 Indeed, we have that, for $y \approx y_0^-$,
  \beq
  \label{95NN}
   \mbox{$
  f''(y)=-2C+ \frac 18\, y_0(y_0-y)+...=0 \,\,\, \mbox{at}
  \,\,\, y_1 =y_0 \big(1 - \frac {16 C}{y_0^2}+...\big), \quad \frac{C}{y_0^2}\ll
  1.
   $}
   \eeq
The second condition in (\ref{94}) then demands
 $$
 \mbox{$
 f(y_1(y_0))=1- \frac {512}3 \, \frac{C^3}{y_0^2}+...>0 \quad
 \big(\frac {C^3}{y_0^2} \ll 1\big).
  $}
  $$
  Existence of $f_1(y)$ is proved if, for some $(y_0,C)$,
 \beq
 \label{94NN}
 y_1=0.
  \eeq

(ii)  We next use again the easy fact that follows from
 the ODE (\ref{4}) that consists of two terms only with clear depending sign
 properties. Thus,  (\ref{85})  implies the positivity,
  \beq
  \label{85N}
  \mbox{$
  -((1-f^2)f''')(y) = \frac 14\, \int_y^{y_0} f' y >0,
   $}
   \eeq
  in a neighbourhood of $y_0^-$,  where $f(y)$ is strictly increasing.
Therefore, $f'''<0$ there, so that $f''(y)$ is strictly decreasing
while $f'(y) \ge 0$.  Hence, we conclude that a monotone function
$f(y)$ cannot have more than a unique\footnote{It seems that this
is
 key for the uniqueness, but we still cannot
 complete the argument, where, as usual,  an extra {\em monotonicity}-like
 property for the ODE is necessary.}
 inflection
point $y_1=y_0(y_0,C)$, at which (\ref{94}) holds, in the
connected domain of monotonicity of $f(y)$.

  (iii) It is easy to
  see from (\ref{k1}) that, for any fixed $y_0>0$,
   \beq
   \label{93}
   f'' <0 \,\,\, \mbox{on} \,\,\, [0,y_0) \,\,\, \mbox{provided that}
   \,\,\, C \gg y_0^2,
    \eeq
so that the inflection point in (\ref{94}) disappears for
sufficiently large $C$.
 It is also easy to see that the type of disappearance of the
 point
 $I(y_0,C)=(y_1(y_0,C),f(y_1(y_0,C)))$ from the domain $S_1=\{y>0, \,\, f \in (0,1)\}$ is
 different for $y_0 \gg 1$ and $y_0 \ll 1$. Namely, for
$y_0 \gg 1$, the point $I(y_0)$ passes through the boundary axis
$\{f=0\}$ as $C$ increases, while for $y_0 \ll 1$, this is done
through the boundary segment on $\{y=0\}$.




 Using the above estimates,
by continuous dependence on parameters $y_0$ and $C$ of this
analytic family $\{f\}$, we conclude that there exists a  $y_0>0$
such that (\ref{94NN}) holds
 and the $f_1(y)$ satisfies (\ref{95}).
$\qed$

\smallskip

 We expect that, in similar lines, there exists a proof of
 existence and uniqueness for higher-order TFEs leading to
 ODEs such as (\ref{4m})
 that still contain two
operators with related ``signs".
 Note that the total
number of parameters increases with the order $2m$ so the shooting
approach becomes more and more involved. For instance, for $m=3$,
i.e., for the sixth-order TFE (\ref{8}), we then have a
three-parameter shooting, with parameters $f'(0)$, $f'''(0)$, and
$f^{(5)}(0)$.
Nevertheless, the uniqueness of $f_1$ is surely observed in
numerical experiments in Section \ref{SectHi}.

\subsection{On stability of the similarity solution in PDE
setting}

We now show how to check whether the similarity solution (\ref{3})
with $f=f_1(y)$ is stable for the TFE (\ref{eq1}) with respect to
small perturbations of the data.  We introduce the rescaled
variables
 \beq
 \label{5.1}
  \mbox{$
u(x,t)= v(y,\t), \quad y=\frac x{t^{1/4}}, \quad \t= \ln t,
 $}
 \eeq
 where $v$ solves the equation with the operator ${\bf B}$ from
 (\ref{4}),
 \beq
 \label{5.2}
 v_\t= {\bf B}(v).
  \eeq
Since $\t \to -\infty$ as $t \to 0$, where $u(x,t)$ approaches a
slightly perturbed data $S_+(x)$, we need to show that these
perturbations are exponentially decaying as $\t$ increases. Thus,
one needs to check the spectrum of the linearized operator
 \beq
 \label{5.3}
  \mbox{$
 {\bf B}'(f)Y= -((1-f^2) Y''')'+ \big(\frac y4 + 2ff'''\big) Y' +
 2(ff''')'Y,
  $}
  \eeq
   that is posed on $y \in (0,y_0)$ with the anti-symmetry regular
   conditions at $y=0$,
 \beq
 \label{5.31}
 Y(0)=Y''(0)=0
  \eeq
    and the free-boundary conditions
 \beq
 \label{5.4}
 Y=Y'=0 \quad \mbox{at} \quad y=y_0.
  \eeq
  For simplicity of linear stability analysis, we assume that the
  free boundary is fixed at $=y_0$, so we exclude
  variation of the free boundary that are not principal.

 Operator (\ref{5.3}) is not symmetric in any space such as
$L^2_\rho$, so it does not admit a self-adjoint extension. It is
rather difficult degenerate singular operator that
 at the singular point $y=y_0$ has the principle part with a
 quadratic degeneracy according to (\ref{k1}),
  $$
  {\bf B}'(f)Y=-2C((y_0-y)^2 Y''')'+... \, .
  $$
Therefore, solving the principle part of the eigenfunction
 equation,
  \beq
  \label{dd1}
 -2C((y_0-y)^2 Y''')'+...= \l Y \approx 0,
  \eeq
 one obtains that the strongest singularity  admitted
 by (\ref{dd1}) is bounded,
  \beq
  \label{dd2}
   Y_{\rm s}(y) \sim (y_0-y) \ln(y_0-y)+... \quad \mbox{for}
   \quad y \approx y_0^-.
 \eeq
Obviously, the singular part (\ref{dd2}) does not satisfy
(\ref{5.4}).

Hence, the conditions (\ref{5.4}) themselves define an operator
extension with a discrete spectrum. This follows from classic
theory of ordinary differential operators \cite{Nai1}, and also
can be easily seen geometrically. Indeed, looking for
eigenfunctions solving (\ref{dd1}) is equivalent to shooting from
the point $y=0$, where in view of (\ref{5.4}) we are left with two
parameters $Y'(0)$ and $Y'''(0)$, and actually, since by  linear
scaling we can always put $Y'(0)=1$ (if it is not zero, which is a
special non-generic case). Therefore, the  parameter $Y'''(0) \in
{\mathbb C}$ together with $\l \in {\mathbb C}$ are, in general,
sufficient to satisfy precisely two conditions (\ref{5.4}) at
$y=y_0$. This gives an analytic system for $\l$ \cite{Nai1}, and
shows that a compact resolvent in a suitable functional space
 as a meromorphic function is  detected.

We restrict our attention to a necessary bound of the spectrum of
(\ref{5.3}),
 \beq
 \label{dd3}
  {\bf B}'(f)Y = \l Y.
   \eeq
 As usual,   multiplying this by the complex conjugate
 ${\overline Y}$ in $L^2$, taking the
 conjugate and multiplying by $Y$, and summing up both yields
  \beq
  \label{dd4}
   \mbox{$
 - \int(1-f^2)\,|Y''|^2 + \int [-(ff')']\,|Y'|^2 + \int
 \big[(ff''')'- \frac 18\big] \, |Y|^2
= \frac {\l+\bar \l}2 \int |Y|^2.
 $}
  \eeq
It is crucial that the first, highest-order term is negative, so
the rest of the  positive terms are expected to be estimated via
interpolation. Actually, there are only two positive components:
(i) in the second term:
 \beq
 \label{dd5}
-(ff')'= - f f'' -(f')^2, \quad \mbox{where} \quad -f f''>0 \,\,\,
\mbox{by} \,\,\, (\ref{95}),
 \eeq
 and (ii) in the last term
\beq
 \label{dd6}
(ff''')'= ff^{(4)}+ f'f''', \quad \mbox{where $f'f'''<0$ but the
sign of $ff^{(4)}$ is unknown}.
 \eeq
Note that close to $y=y_0$, the interpolation is possible
regardless the degeneracy of the principle part in (\ref{dd4}).

Thus we conclude from (\ref{dd4}) that
 \beq
 \label{dd7}
 {\rm Re}\, \l < 0,
  \eeq
  provided that  a suitable (Friedrichs) self-adjoint extension of
the following symmetric operator satisfies (a proper setting is
standard, \cite{Nai1}):
 \beq
 \label{dd8}
  \mbox{$
 {\bf Q}_P(f)Y= - [(1-f^2)Y'']''+ [(f f')'Y']' + \big[(f f''')'-
 \frac 18 \big]Y <0.
  $}
  \eeq
For the first FBP profile $f_1(y)$, which is not given explicitly,
(\ref{dd8}) can be checked numerically and is rather plausible. At
least, even without careful and not that easy spectral numerics,
we observe that the self-adjoint operator in (\ref{dd8}) exhibits
a clear tendency to be  negatively determined for reasonable FBP
profiles $f=f_1(y)$ from (\ref{95}).

\section{The Cauchy setting for RP--$1$:  oscillations and  similarity profiles}
 \label{Sect31}

We begin by noting that the CP setting for the TFE (\ref{1}) {\em
assumes no analyticity}
 unlike the TFE setting discussed in
Section \ref{Sect2}. For the CP, we observe oscillatory and
changing sign behaviour near the interfaces.

\subsection{Local solutions of changing sign}

The oscillatory local behaviour for (\ref{4}) is more complicated.
Here we follow \cite{Gl4}. We introduce  the oscillatory component
$\varphi(s)$ as:
 \beq
 \label{k3}
 f(y) = 1-(y_0-y)^3 \varphi(s), \quad s= \ln (y_0-y),
  \eeq
  where $\varphi(s)$ solves an autonomous ODE with a discontinuous
  nonlinearity,
   \beq
 \label{k4}
  \textstyle{
P_3(\varphi) \equiv  \varphi''' + 6 \varphi'' + 11 \varphi' +
6\varphi =-\l_0 \,{\rm sign} \, \varphi,
  \quad \mbox{where} \quad
\l_0= \frac 14 \, y_0>0.
 }
 \eeq
 We look for a periodic solution of (\ref{k4}), which
by (\ref{k3}) describes the oscillatory nature of solutions in the
CP. Let us first note that  technique based on calculation of
rotation of vector fields for the dynamical system (\ref{k4}) does
not apply here; see
 \cite[p.~45-53]{KrasZ}. In this case, existence of multiple periodic solutions
 depends on the careful analysis of the index of a guiding
 function $V$ (it always exists  and can be calculated explicitly
 by using the linear form of the operator at infinity).
 Indices of such periodic solutions are unknown, the main results
 in \cite[p.~52]{KrasZ} seem not applicable for present problem.

 Therefore,  a direct shooting approach is effective here.
 We consider a 2D bundle of  orbits $\varphi(s,C_1,C_2)$ satisfying the conditions
  \beq
  \label{v1}
  \varphi(0)=0, \quad  \varphi'(0)=C_1, \quad  \varphi''(0)=C_2,
   \eeq
   where $C_{1,2} \in \re$ are arbitrary parameters.
  Let us state some key properties implying existence of a periodic orbit.

  (i) {\em $\varphi(s,C_1,C_2)$ is uniformly bounded for all $s >0$.} Obviously, in view of
  boundedness of the nonlinearity, $\varphi(s,C_1,C_2)$ is locally well
 defined. If, on the contrary, $\varphi(s,C_1,C_2)$ is unbounded,
  we should have that the linear counterpart,
  $$
   \varphi''' + 6 \varphi'' + 11 \varphi' +
6\varphi = 0, $$
 must admit unbounded solutions. Setting
$\varphi(s) = {\mathrm e}^{\mu s}$ gives  the characteristic
equation
 $$
 \Phi(\mu) \equiv \mu^3 + 6 \mu^2 + 11 \mu + 6=0.
 $$
 Since, $\Phi'(\mu)>0$ for $\mu \ge 0$, $\Phi(0)=6$,
 and
  $$
  \mbox{$
  \Phi'(\mu) = 0 \quad \mbox{for} \quad \mu=\mu_\pm= \frac 13(6
  \pm \sqrt 3), \quad
  \mbox{with} \quad
  \Phi(\mu_\pm) = \mp 0.3849... \, ,
   $}
   $$
  we have that all three roots of $\Phi(\mu)$ are negative.
   So the solutions cannot grow without bound as $s \to + \infty$.

   Actually, this means that, as a consequence, that the fourth-order ODE (\ref{k4})
   with bounded coefficients is a dissipative DS
having a bounded absorbing set.
 Dissipative DSs are known to admit periodic solutions in a rather
general setting \cite[\S~39]{KrasZ} provided these are
non-autonomous (so the  period is fixed). For the autonomous
system (\ref{k4}), the proof in \cite[\S~7.1]{Gl4} can be
completed by  shooting.

 (ii) {\em  $\varphi(s,C_1,C_2)$ is oscillatory as $s \to +\infty$, i.e., has
infinite number of zeros in any neighbourhood of $s= +\infty$}.
Indeed,  if, say, $\varphi(s,C_1,C_2)>0$
 for all $s \gg 1$, we will have there the linear ODE
  $$
   \varphi''' + 6 \varphi'' + 11 \varphi' +
6\varphi = -1,
$$
admitting the particular equilibrium solution $\varphi = - \frac 16$. Hence, by the above stability
of the orbits
of the homogeneous linear equation in (i), we conclude that  $\varphi(s,C_1,C_2) \approx
- \frac 16 <0$ for all $s \gg 1$, from whence comes a contradiction.

We now state the main result.


 \begin{proposition}
  \label{PrPer}
  {\rm (i)} Equation $(\ref{k4})$ admits a unique non-trivial periodic solution
  $\varphi_0(s)$ of changing sign, which is  stable as $s \to + \infty$; and

  \noi {\rm (ii)} The periodic orbit $\varphi_0(s)$ of
  $(\ref{k4})$ is the unique bounded connection with $s= -\infty$.
   \end{proposition}

  \noi {\em Proof.} (i) Existence and uniqueness have a pure algebraic proof,
\cite[\S~7.4]{Gl4}. Stability follows from the above analysis of
the linearized problem.  (ii)  This also implies that the stable
manifold of $\varphi_0(s)$ as $s \to -\infty$ is empty, so, by the
geometric analysis in \cite{Gl4}, $\varphi_0(s)$ is the only
bounded
  orbit connecting $s=-\infty$. $\qed$

  \smallskip

 Figure \ref{F16}  shows that the periodic solution
 $\varphi_0(s)$ is  {\em stable},
 up to translations in $s$, as usual, for the
 ODE (\ref{k4}) as $s \to + \infty$.




 In view of (ii), the expansion (\ref{k3}) with  $\varphi=\varphi_0(s+s_0)$, $s_0 \in \re$ is arbitrary,
 is the only maximal regularity connection with $f(y) \equiv 0$ for $y>y_0$.
 Therefore, this describes the generic robust structure of the multiple zero at the interface
 of arbitrary similarity solutions of the TFE--4 under consideration.
 We expect that this structure remains unchanged for general solutions of the PDE
  (\ref{1}), which is a difficult open problem.
  Observe that according to (\ref{k3}), the actual
 regularity of $f(y)$ is $C^{2,1}$ that is better than for the
 above FBP with the less smooth $C^{1,1}$ behaviour in (\ref{k1})
 (if we formally set $f(y) \equiv 1$ for $y \ge y_0$).
As we have mentioned, the regularity in (\ref{k3}) can be
attributed to the Cauchy problems and in fact is the {\em maximal
regularity} which is admitted by the ODE under consideration,
\cite{Gl4}.


\begin{figure}
\centering
\includegraphics[scale=0.60]{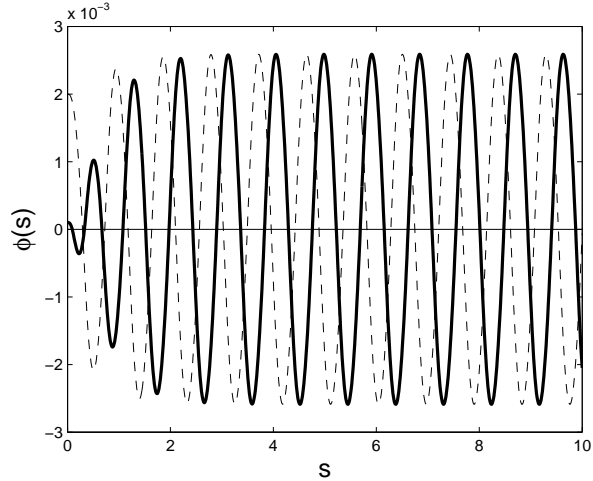}  
\vskip -.3cm \caption{\small On the stable periodic orbit of
(\ref{k4}), with $y_0=4$.}
   \vskip -.3cm
 \label{F16}
\end{figure}


\subsection{Similarity profile $f_0(y)$ for the CP}

We begin by noting that finite propagation in the ODEs such as
(\ref{4}), i.e., existence of a finite interface $y_0$, is
well-known for a long time; see first ODE proofs in \cite{Bern88,
BMc91}, and more general energy estimates for PDEs in
\cite{Bern01, Shi2} and survey in \cite{GS1S-V}.

As in the FBP, the bundle of asymptotics (\ref{k3}) is 2D
comprising parameter $\{y_0,s_0\}$, where $s_0$ is the translational
constant  in the oscillatory component, $\varphi_0(s) \mapsto
\varphi_0(s+s_0)$. Of course, the matching problem for the CP with the
oscillatory bundle (\ref{k3}) gets more difficult.

In Figure \ref{Fch1}, we show the unique similarity profile
$f_0(y)$ corresponding to the CP. Figure \ref{Fch2} explains the
character of oscillations of the $f(y)$ about $f= 1$. This
behaviour well-corresponds to the expansion (\ref{k3}). This shows
a rough estimate of the interface location
 \beq
   \label{ii11}
   y_{0,{\rm CP}} \sim 6.7,
    \eeq
which is difficult to improve numerically in view of the
oscillatory behaviour, the necessity of regularization parameters,
and sufficiently low tolerance of convergence that was applied to
the  dynamical system (equivalent to (\ref{4})) with the
essentially non-symmetric matrix.

\begin{figure}
\centering
\includegraphics[scale=0.65]{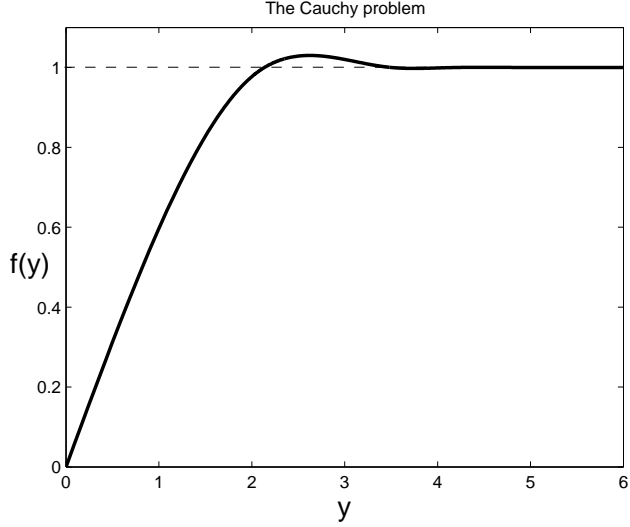}  
\vskip -.3cm \caption{\small The  CP profile $f_0(y)$ satisfying
(\ref{4}), (\ref{7}).}
   \vskip -.3cm
 \label{Fch1}
\end{figure}


\begin{figure}
\centering
\subfigure[oscillations, enlarged]{
\includegraphics[scale=0.52]{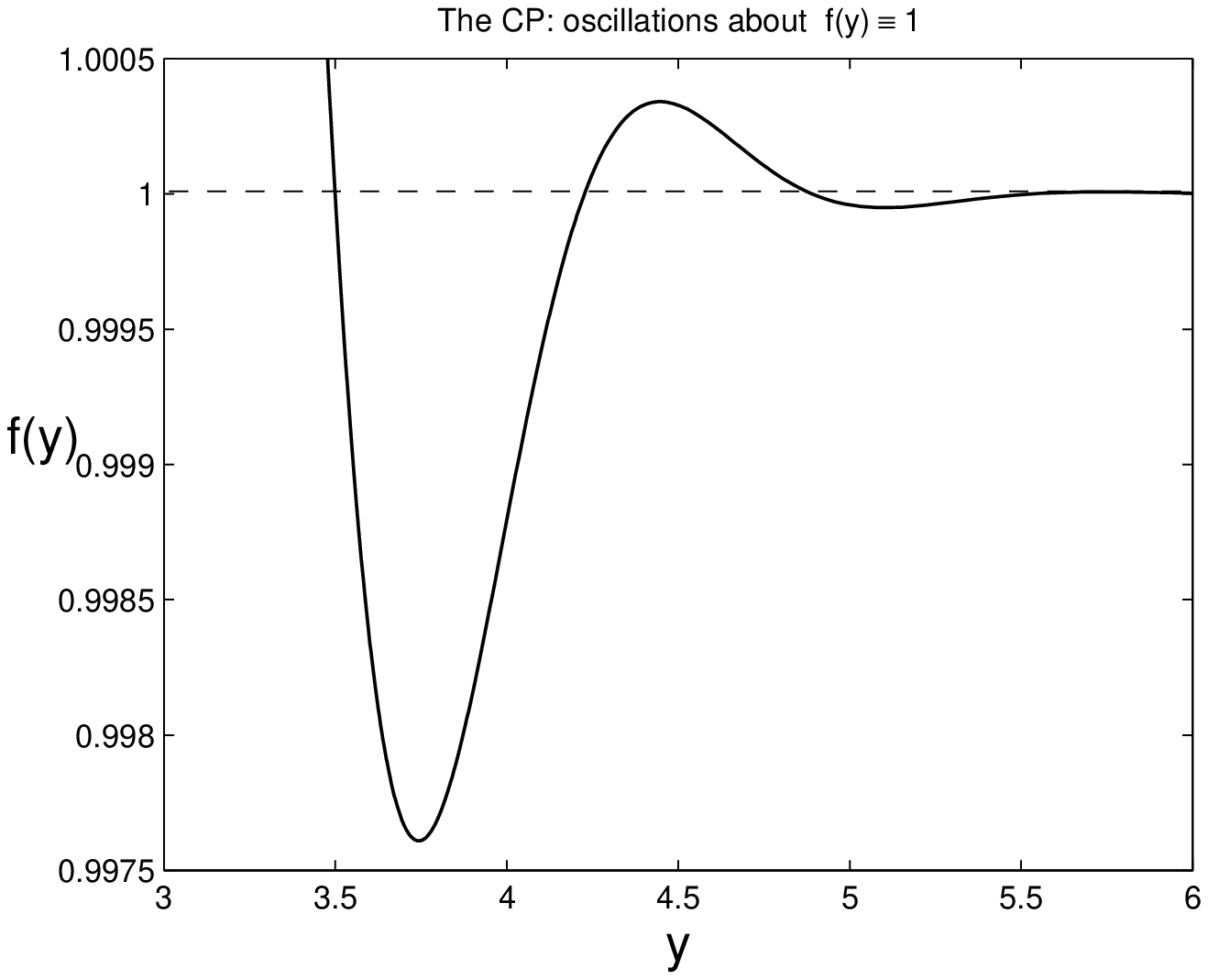} 
}
\subfigure[oscillations, further enlarged]{
\includegraphics[scale=0.52]{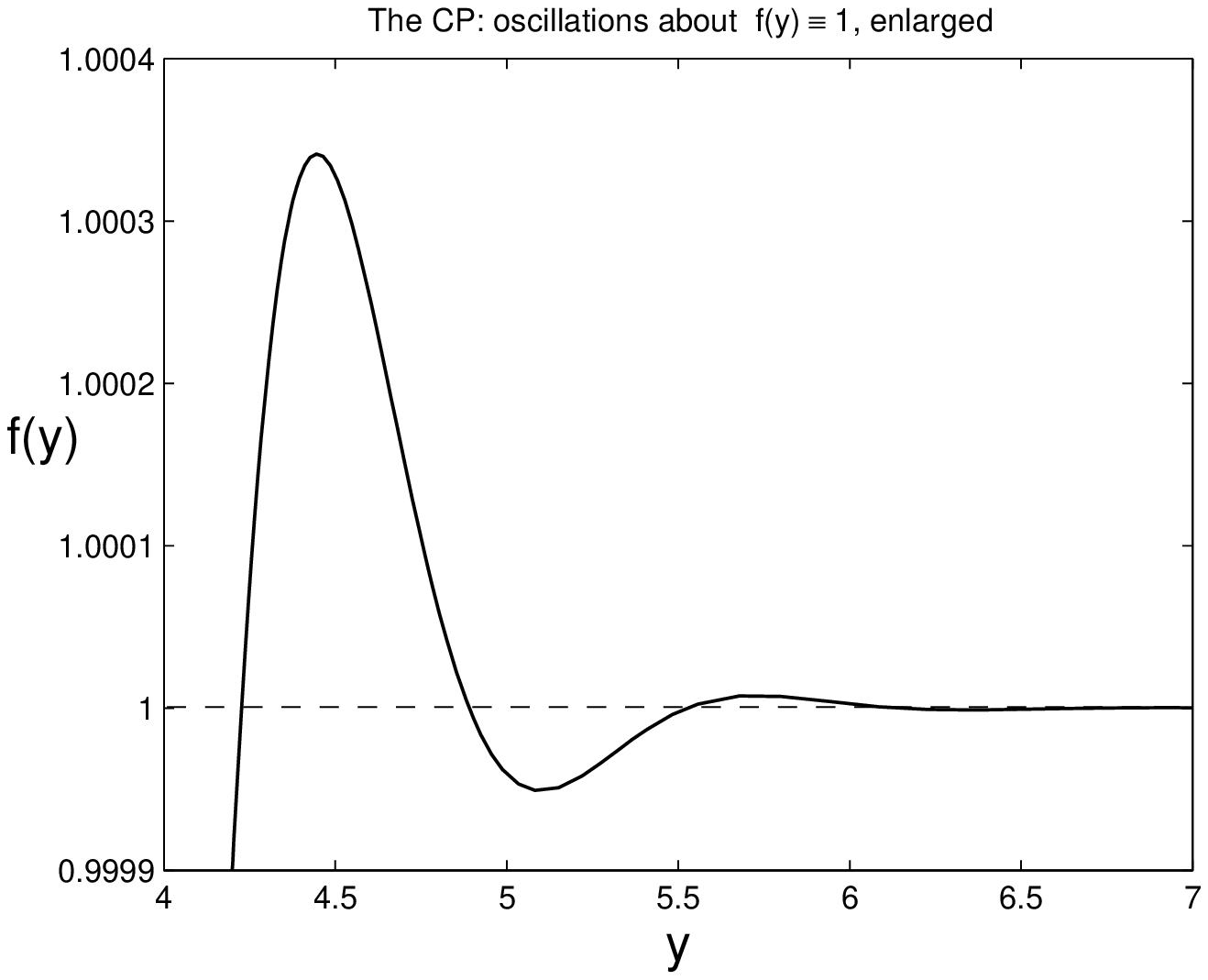} 
}
   \vskip -.3cm
    \caption{\small  Enlarged oscillations of $f_0(y)$ from  Figure
\ref{Fch1} about the equilibrium $f=1$.}
   \vskip -.3cm
 \label{Fch2}
\end{figure}


\subsection{On a countable set of FBP profiles}

Figure \ref{F01} indicates an interesting ``geometric" relation
between the first FBP profile and the CP one ($n=1$).
 For further convenience, in Figure \ref{F02N},
 we present the enlarged version of this figure with extra
 details, where we include shooting of the third FBP profile $f_3(y)$
 with interface positions at, respectively,
  $$
  y_{01} \sim 2.35, \quad y_{02} \sim 3.75, \quad \mbox{and} \quad
  y_{03} \sim 4.74,
  $$
  which all are smaller than the CP interface location
  (\ref{ii11}).
By boldface dotted line therein, we denote the CP profile
$f_0(y)$.

\begin{figure}
\centering
\includegraphics[scale=0.75]{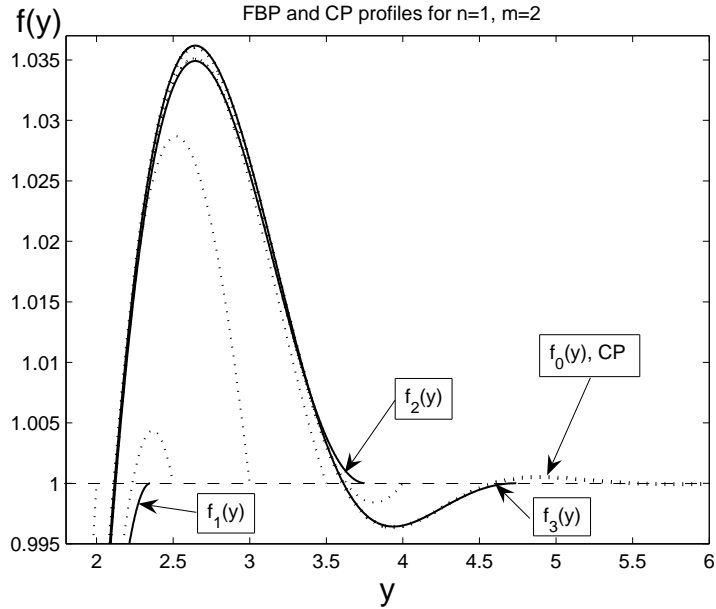}  
\vskip -.3cm \caption{\small The  CP  $f_0(y)$ and first three FBP
profiles for $n=1$.}
   \vskip -.3cm
 \label{F02N}
\end{figure}

 Namely, we see that the touching point
of $f_1(y)$ at the level $1$ belongs to the interval of the first
oscillation of the profile for the CP, and the same is true for
$f_2$ and $f_3$ in Figure \ref{F02N}. Moreover, we expect a
further extension of this geometric property (an open problem).

\smallskip

\noi{\bf Conjecture \ref{Sect31}.2.} {\em The FBP
$(\ref{4})$--$(\ref{6})$ admits a countable family of different
solutions $\{f_k(y)\}$, where $f_1(y) \in [-1,1]$, which converges
to the CP profile:
 \beq
 \label{CP1}
  f_k(y) \to f_0(y) \quad \mbox{as \,\, $k \to \infty$}
   \eeq
uniformly on any interval from $[0,y_{0,{\rm CP}})$.}

\smallskip

In Figure \ref{Fcount}, for clarity,  we  present a formal scheme
of this geometric interaction between FBP and CP profiles, which
partially has been reflected in Figure \ref{F02N}. This was
already seen in Figures \ref{F01}, \ref{F02N},
 and will be clearer explained in Figure
\ref{Figm39} obtained numerically for the TFE-6 ($m=3$ in
(\ref{8})) admitting more oscillatory patterns. In particular,
these figures show that the interface points $y_{0k}$ of the FBP
profiles satisfy
 \beq
 \label{CP1k}
 y_{0k} \to y_{0,{\rm CP}}^- \quad \mbox{as} \quad k \to \infty.
  \eeq

  Recall that any FBP profile $f_k(y)$ being put into the similarity formula (\ref{3})
  gives a solution of the problem (\ref{1}), (\ref{2}).
Therefore, in the FBP setting, {\em Riemann's problem for the
TFE--4 has an infinite countable number of solutions} (and only
the first profile $f_1(y)$ has the physical range $[-1,1]$), while
for the Cauchy setting, {\em such a solution is expected to be
unique} (and violates the range $[-1,1]$).

\begin{figure}
 \centering
 \psfrag{f(y)}{$f(y)$}
 \psfrag{f_1(y), FBP}{$f_1(y)$, FBP}
 \psfrag{f_2(y), FBP}{$f_2(y)$, FBP}
 \psfrag{f_3(y), FBP}{$f_3(y)$, FBP}
 \psfrag{f_4(y), FBP}{$f_4(y)$, FBP}
 \psfrag{f(y), CP}{$f(y)$, CP}
 \psfrag{y_0CP}{$y_{0,{\rm CP}}$}
 \psfrag{1}{$1$}
 \psfrag{t4}{$t_4$}
  \psfrag{v(x,t-)}{$v(x,T^-)$}
  \psfrag{final-time profile}{final-time profile}
   \psfrag{tapp1}{$t \approx 1^-$}
\psfrag{y}{$y$}
 \psfrag{0<t1<t2<t3<t4}{$0<t_1<t_2<t_3<t_4$}
  \psfrag{0}{$0$}
 \psfrag{l}{$l$}
 \psfrag{-l}{$-l$}
\includegraphics[scale=0.43]{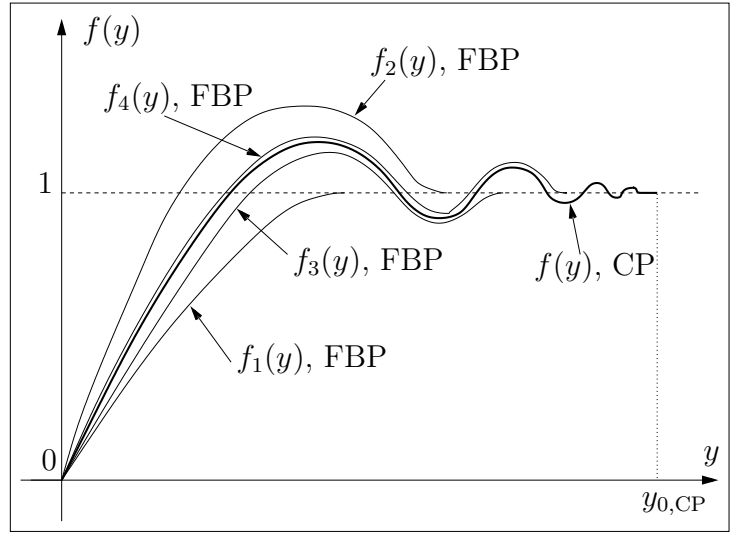}     
\caption{\small An extra schematic explanation of a countable set
$\{f_k(y)\}$ of FBP profiles approaching the CP one (the bold
line); cf. Figure \ref{F02N}.}
     \vskip -.3cm
 \label{Fcount}
\end{figure}

\section{On extension at the branching point $n=0$ in the Cauchy problem}
 \label{SectBr}

Here, we demonstrate another approach to Riemann's problem in the
case of the Cauchy problem for the TFE-4 (\ref{n1}) with parameter
$n \ge 0$. Then the  similarity solution has the same form
(\ref{3}), there $f(y)$ solves the ODE
 \beq
 \label{s1}
  \begin{matrix}
 {\bf B}_n(f) \equiv -(|1-f^2|^n f''')' + \frac 14 \, f'y =0
 \quad \mbox{for} \quad y>0,\smallskip\smallskip\smallskip \\
 f(0)=f''(0)=0, \quad f(+\infty)=1.
  \end{matrix}
   \eeq
By our local analysis, we expect the solution with the finite
interface to have a sufficiently regular connection with the
equilibrium $f = 1$, as explained in Section \ref{Sect31}.

Our basic idea is to show that the solution of (\ref{s1}) can be
extended from the solution $f_0$ in (\ref{b3}) of the linear
problem for $n=0$, which then becomes a {\em branching point} for
(\ref{s1}); see \cite[p.~371]{KrasZ}. To justify this,  we will
need some extra spectral properties of the linear operators
involved.

 \subsection{Similarity profile and spectral properties of ${\bf B}$
 for $n=0$}
 \label{SectExi}

 Thus, as we have seen, for $n=0$, the problem (\ref{s1}) in the
 CP setting takes the form
\beq
 \label{s10}
  \begin{matrix}
 {\bf B}_0 (f) \equiv -f^{(4)} + \frac 14 \, f'y =0
 \quad \mbox{for} \quad y>0, \smallskip\smallskip\ssk \\
 f(0)=f''(0)=0, \quad f(+\infty)=1,
  \end{matrix}
   \eeq
and has
  the unique solution $f_0(y)$ given in (\ref{b3}).
It follows from (\ref{s10}) that
 \beq
 \label{bb1}
 \mbox{$
 {\bf B}_0= {\bf B}- \frac 14 \, I,
  $}
  \eeq
  where ${\bf B}$ is the linear operator in (\ref{b2}) that
  defines the rescaled kernel of the fundamental solution of the parabolic operator
   $D_t
  + D_x^4$.

   The necessary spectral properties of
the linear non self-adjoint operator ${\bf B}$ introduced in
(\ref{b2}) and the corresponding adjoint operator ${\bf B}^*$ are
of importance in the  asymptotic analysis, as
  explained in \cite{Eg4} for similar general $2m$th-order
 operators (see also \cite[\S~4]{Bl4}).
 In particular, ${\bf B}$ is naturally defined
in a weighted space $L^2_\rho(\re)$, with
 $$
 \rho(y)={\mathrm e}^{a
|y|^{4/3}}, \quad \mbox{where \,$a \in (0,2d)$ \, is a constant},
 $$
 with the domain being the corresponding Hilbert (Sobolev) space
 $H^4_\rho(\re)$. Then ${\bf B}: H^4_\rho(\re) \to L^2_\rho(\re)$
 is a bounded operator and
    has the discrete (point) spectrum
 \begin{equation}
\label{spec1}
  \sigma({\mathbf B}) = \big\{\lambda_l = -\mbox{$\frac
l4$}, \,\, l = 0,1,2,...\big\}.
\end{equation}
The corresponding eigenfunctions are normalized derivatives of the
rescaled kernel,
 \beq
 \label{psi1}
 \mbox{$
 \psi_l(y)= \frac {(-1)^l}{\sqrt{l!}} F^{(l)}(y), \quad l=0,1,2,... \,
 .
  $}
  \eeq

The adjoint operator
 \beq
 \label{ad1}
  \mbox{$
 {\bf B}^*=- D_y^4 - \frac 14 \, y D_y
  $}
  \eeq
  has the same spectrum (\ref{spec1}) and the polynomial
  eigenfunctions
  \begin{equation}
 \label{psi**1}
  \mbox{$
 \psi_l^*(y) = \frac 1{\sqrt{l !}}
\sum_{j=0}^{\lfloor-\lambda_l\rfloor} \frac {1}{j !}D^{4j}_y y^l,
\quad l=0,1,2,...,
 $}
 \end{equation}
which form a complete subset  in $L^2_{\rho^*}(\re)$, where
$\rho^*=\frac 1\rho$. Similarly, the domain of the bounded
operator ${\bf B}^*$ is the Sobolev space $H^4_{\rho^*}(\re)$. In
particular, it is easy to compute
   \beq
   \label{BB3}
    \mbox{$
 \psi_0^*=1, \,\, \psi_1^*=  y, \,\,\,
 \psi_2^*=\frac 1{\sqrt 2} \, y^2, \,\,\,
 \psi_3^*=\frac 1{\sqrt 6} \, y^3, \,\,\,
   \psi^*_4= \frac 1{\sqrt{24}} \,(y^4+24),
 $}
    \eeq
    etc.
 As ${\bf B}$, the adjoint operator ${\bf B}^*$
has the compact resolvent $({\bf B}^*-\lambda I)^{-1}$ in
$L^2_{\rho^*}(\re)$. It is not difficult to see that, on
integration by parts, the eigenfunctions
 (\ref{psi1})
 are orthonormal to polynomial eigenfunctions
$\{\psi^*_l\}$ of the adjoint operator ${\bf B}^*$, so
 \beq
 \label{ort1}
 \langle \psi_l, \psi^*_k \rangle = \d_{lk},
  \eeq
  where $\langle \cdot, \cdot \rangle$ denotes the standard (dual)
  scalar product in $L^2(\re)$.

\smallskip

Thus, in the necessary restriction  to odd functions denoted by
$L^2_\rho(\re_+)$, according to (\ref{bb1}), we have that ${\bf
B}_0$ has the discrete spectrum
 \beq
 \label{ss1}
  \mbox{$
  \s({\bf B}_0)= \big\{- \frac {1+l}4 \le - \frac 12, \,\, l=1,3,5,... \big\}
   $}
 \eeq
 and a complete and closed set of eigenfunctions $\{\psi_l\}$.
 Hence,  the linearized operator at $n=0$ is strictly negative in
 view of (\ref{ss1}) that gives a good chance to extend the
 solution $f_0$ for $n>0$ small enough by applying classic
 branching theory
 in the study of the limit $n \to 0^+$ in the problem (\ref{s1});
see \cite[p.~319]{VainbergTr} and \cite[p.~381]{Deim}.

 \subsection{Branching at $n=0$}

To this end,
    for $n>0$, we write down (\ref{s1})
 as a perturbation of the linear problem (\ref{s10}),
 \beq
 \label{ne}
  \mbox{$
 {\bf B}_0 f = g(f;n) \equiv  \big[(|1-f^2|^n-1)f''' \big]'.
  $}
  \eeq
Next, since ${\bf B}_0<0$ is known to have compact resolvent in
$L^2_\rho(\re_+)$, we  consider the equivalent integral equation
for the function
 \beq
 \label{nn1NN}
 Y=f-f_0, \quad \mbox{i.e.,} \quad f= f_0+Y,
  \eeq
 which has the form
  \beq
  \label{ne1}
  Y={\bf A}(f_0+Y;n) \equiv {\bf B}_0^{-1} g(f_0+Y;n).
  \eeq
  Initially, we formally assume that,
 for small $n \le 1$, the nonlinear operator on the right-hand
side can be treated as a compact Hammerstein operator in
$L^p_\rho$-spaces, as classic theory suggests \cite[\S~17]{KrasZ}.
Recall that for $n=0$, (\ref{ne1}) is indeed a linear integral
equation with a compact operator admitting the unique (up to a
multiplier) solution $f_0$ as in (\ref{b3}).


 Therefore, performing below necessary
computations, we bear in mind a formal use of classic branching
theory for compact integral operators  as in (\ref{ne1}). We also
do not discuss here specific aspects of compact integral operators
in weighted $L^p$-spaces; see \cite{Eg4} for extra details.
 Note that, in a suitable metric, ${\bf B}$ is a sectorial
 operator, \cite[\S~5.1]{2mSturm}.

  Therefore,
for derivation of branching equations, we  use general branching
theory in Banach spaces; see Vainberg--Trenogin
\cite[Ch.~7]{VainbergTr}. Beforehand, we need to discuss a few
typical difficulties associated with the above eigenvalue problem.

As (\ref{ne}) suggests,
 the crucial part is played by the nonlinearity
  \beq
  \label{nn1}
 Q(f;n)=|1-f^2|^n-1 \quad \mbox{in a neighbourhood of the point
 $P_0=\{f=1,n=0^+\}$}.
  \eeq
One can see that the derivative $Q'_f(f;n)$ {\em is not}
continuous at
 $P_0$.
 Therefore, the standard
  implicit function (operator) theorem
  \cite[p.~319]{VainbergTr} does not apply.

Nevertheless, the application of index-degree theory
\cite[p.~355]{KrasZ} formally  demands the differentiability only
at the given point $f=f_0$, $n=0$, which is true (cf. more
delicate computations below).
 Indeed, for $n=0$,
  $g(f,0) \equiv 0$, and the perturbation disappears in
  (\ref{ne1}) at the branching point.

On the other hand,  in view of our difficulties with the
regularity (and also with
 compactness of the operators involved), it is better to rely on Theorem 28.1
 in \cite[p.~381]{Deim} that is formulated in the linearized
 setting, where the differentiability is ``replaced" by the control of
 higher-order nonlinear terms as in (\ref{nn1}) in a neighborhood.
  As usual, the key principle of branching
 is that
   the corresponding
 eigenvalue has odd multiplicity that can be easily checked in some cases
 (the even multiplicity case needs an additional treatment, which is also a
   routine procedure not to be treated here);
 see further comments below.
 Notice that the condition on the nonlinearity (\ref{nn1}) is
 rather tricky to check that it satisfies the estimate (b) in Theorem 28.1
 in \cite[p.~381]{Deim} for the variable $x$ in $f= F+x$.

 Thus, we can proceed and conclude that
(\ref{s1}) admits a solution $f(y)$ for all sufficiently small
$n>0$ and moreover these profiles form a continuous curve (an
$n$-branch).

We next detect a precise behaviour of this branch as $n \to 0^+$.


 \subsection{
 Asymptotic expansion of the branch  for small $n>0$}

  We need to use in (\ref{ne}) the following expansion:
   \beq
   \label{Ex.6}
   |1-f^2|^n-1= n \ln |1- f^2| + o(n)
 \quad \mbox{as}
 \quad n \to 0^+,
  \eeq
which of course is not true uniformly on bounded intervals in $f$
and should be understood
 in the weak sense. This is crucial for the equivalent integral equation
  (\ref{ne1}).

\begin{proposition}
 \label{Pr.W}

 For the function $f=f_0$ given by $(\ref{s10})$,
 in the sense of distributions and in the weak sense in
$L^\infty(\re_+)$ $($and also in the sense of bounded measures in
 $\re_+)$
 \beq
 \label{Is.1}
 \mbox{$
  \frac 1n \,\big(\, |1-f^2|^n-1\,\big) \rightharpoonup   \ln |1-f^2| \quad \mbox{as} \quad
  n \to 0^+.
   $}
   \eeq
   \end{proposition}

   According to (\ref{Is.1}), analyzing the integral equation (\ref{s10}), we can use
   the fact that,
 for any function $\phi \in {\mathcal L}$, $\phi \in L^1(\re_+)$  (or $\phi \in
 C_0(\re_+)$)
 \beq
 \label{Is.2}
  \mbox{$
  \int (|1-f_0^2|^n-1) \phi(y) \, {\mathrm d}y =
n \big[\int f(y)\ln |1-f_0^2(y)| \,\phi(y) \, {\mathrm d}y + o(1)
\big] \,\,\, \mbox{as} \,\,\, n \to 0^+.
 $}
 \eeq

\smallskip

Let us finish our formal computations that, for convenience, are
performed for the differential equation (\ref{s1}). Namely,
 substituting expansions
(\ref{Ex.6}) into (\ref{s1}) and using (\ref{nn1NN})
 yields the following
perturbed equality:
 \beq
 \label{An.1}
  \mbox{$
 {\bf B}_0 Y= n(\ln |1-f_0^2) f_0 ''')' + o(n),
  $}
  \eeq
and this gives the unique solution
 \beq
 \label{An.22}
 Y= n{\bf B}_0^{-1} \big[(\ln |1-f_0^2| \,  f_0 ''')'  \big] + o(n)
  \quad \mbox{as} \quad n \to 0.
   \eeq
One can see that
 $$
 (\ln |1-f_0^2| \,  f_0 ''')' \in L^2_\rho(\re_+),
 $$
 so (\ref{An.22}) makes sense.
Thus, (\ref{nn1NN}) and (\ref{An.22}) actually determine the
behaviour of the $n$-branch of similarity profiles of Riemann's
problem for sufficiently small $n>0$.


\subsection{On global continuation of $n$-branches: open problem}

Global continuation of branches of nonlinear eigenfunctions of
(\ref{s1}) from the branching point is often an intriguing open
problem. Global bifurcation results concerning continuous branches
of solutions originated at $n=0$ are already given in
Krasnosel'skii (the first Russian edition was published in 1956),
\cite[p.~196]{Kras}. Concerning further results and extensions,
see references in \cite[Ch.~10]{Deim} (especially, see
\cite[p.~401]{Deim} for typical global continuation of bifurcation
branches), and also \cite[\S~56.4]{KrasZ}.

 In general, for the integral equation (\ref{ne1}) with
compact Hammerstein operators in weighted $L^2$-spaces, it is
known since Rabinowitz's study (1971) that  branches are
infinitely extensible and can end up at further bifurcation
points; see \cite[\S~29]{Deim}  for further information stated in
the framework of bifurcation analysis and \cite[\S~56.4]{KrasZ}.
The above $n$-extension of the branch is possible at any $n=n_0>0$
provided that the linearized operator has proper spectral
properties. For oscillatory profiles $f(y)$ about $f=1$, with a
complicated infinite set of intersection points, this is difficult
to check in general but, in principle, can be done in specific
weighted spaces, provided
 assuming that the
oscillatory structure near interfaces is known in detail.
Nevertheless, for the non-variational eigenvalue problem
(\ref{s1}) with non-divergent and non-monotone operators, a
rigorous treatment of global behaviour of branches is very
difficult and remains open.

Since  we have obtained the unique continuous branch originated at
$n=0$, it will never come back to the origin, so this suggests
that it can be extended up to $n=1$ and even further until
homoclinic and other bifurcation points will destroy its necessary
quality. Nevertheless, we do not know principles of such a global
continuation, so we end up our branching analysis as follows:

\smallskip

\noi{\bf Conjecture \ref{SectBr}.1.}  {\em The $n$-branch of
nonlinear eigenfunctions of $(\ref{s1})$ that is originated at
$n=0$ from $f_0$ exists for all $n \in [0,1]$ and  does not have
turning $($saddle-node$)$ points.}

 \smallskip

In Figure \ref{FigBr}, we show the deformation with $n$, with the
step $\D n=0.1$, of the CP similarity profile $f_0(y)$ satisfying
(\ref{s1}).
 It is seen that
all the profiles have quite similar shapes for $n=0$ (the dotted
line) and $n=1$ including also the case of the negative exponents
 $$
  n = -0.1, \,\, -0.2, \,\, -0.3, \,\, -0.4, \,\, -0.5,
  $$
  for which $f_0(y)$ does not have finite interface but remains
  equally oscillatory as $y \to + \infty$.
Figure \ref{FigBr1} shows the enlarged oscillations about the
constant equilibrium 1 of the profiles from Figure \ref{FigBr}.

\begin{figure}
\centering
\includegraphics[scale=0.65]{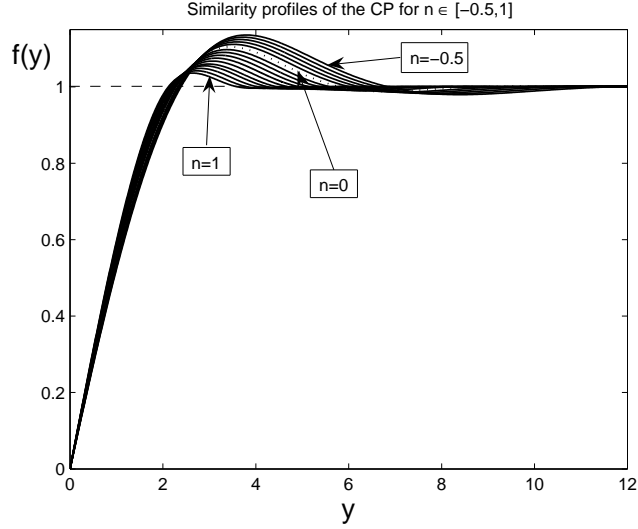}
\vskip -.3cm
 \caption{\small Continuous deformation in $n$
  of  similarity profiles satisfying (\ref{s1}); $ n \in [-\frac 12,1]$.}
   \vskip -.3cm
 \label{FigBr}
\end{figure}

\begin{figure}
\centering
\includegraphics[scale=0.70]{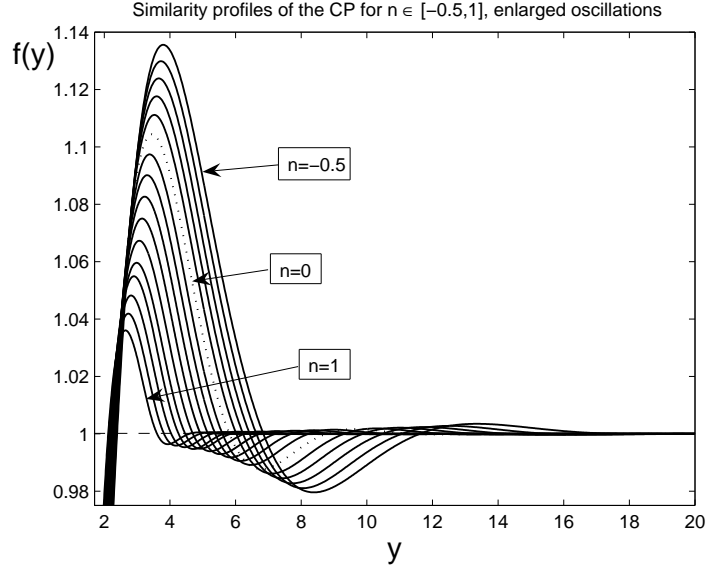}  
\vskip -.3cm
 \caption{\small Enlarged oscillations of    similarity profiles
 from Figure \ref{FigBr};
 $ n \in [-\frac 12,1]$.}
   \vskip -.3cm
 \label{FigBr1}
\end{figure}


\section{Riemann problem $\frac 12$: similarity solutions}
 \label{Sect12}

We now consider the RP--$\frac 12$ (\ref{4}), (\ref{6}),
(\ref{5R}) or (\ref{7}). Key aspects of the similarity analysis
(both theoretical and numerical) are very similar.
  We  stress the
attention to some distinctive features.

In Figure \ref{FCP1}, we show the unique similarity profile
$f_0(y)$ of the RP--$\frac12$ in  the Cauchy setting. Notice that
at the equilibrium level $f=1$ the similarity profile is clearly
less oscillatory than at $\{f=0\}$, where propagation is
approximately governed by the linear bi-harmonic PDE (\ref{l1}).

\begin{figure}
\centering
\includegraphics[scale=0.65]{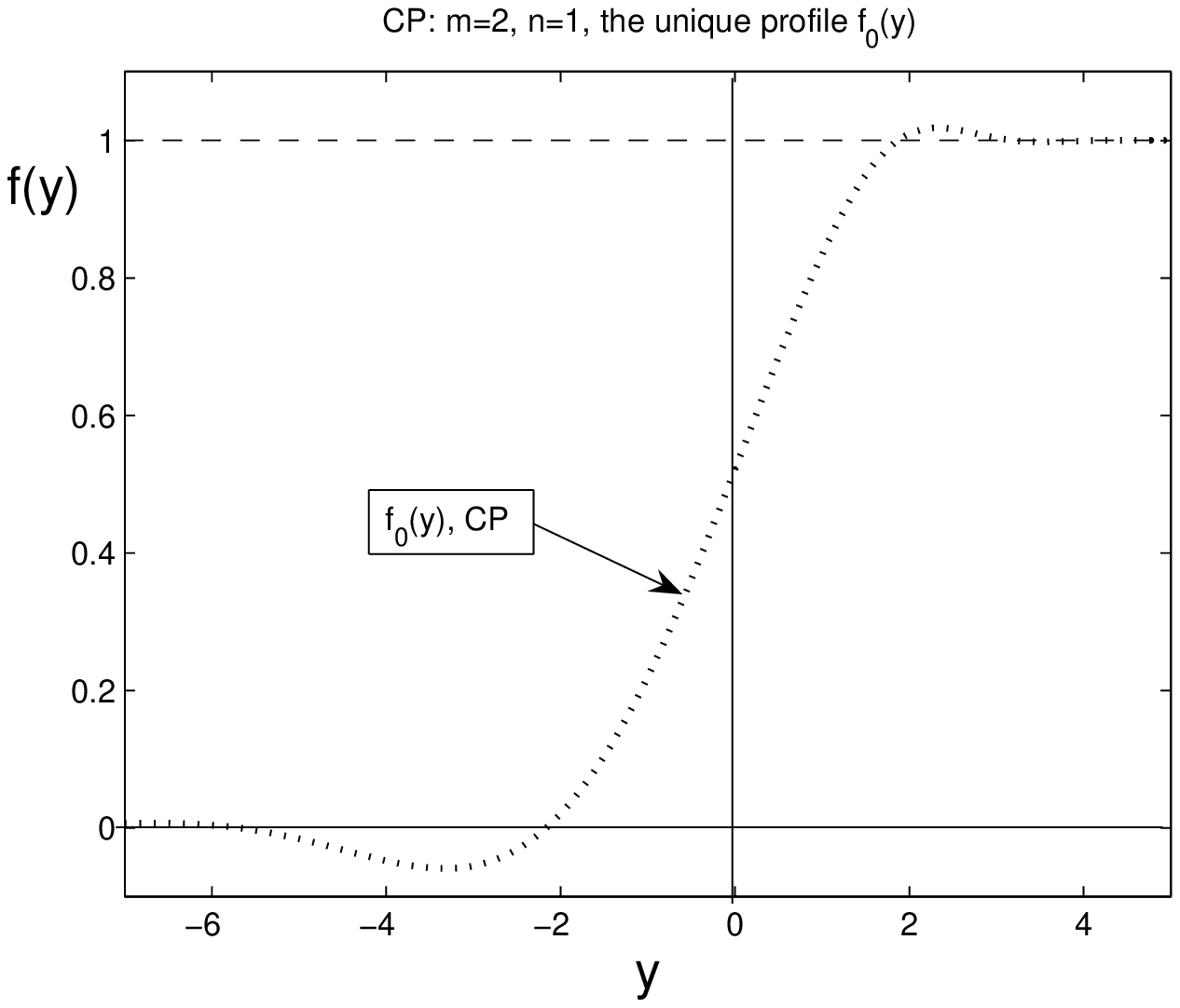}  
\vskip -.3cm
 \caption{\small The CP solution of the RP--$\frac 12$
(\ref{4}), (\ref{6}), (\ref{7}) for $n=1$.}
   \vskip -.3cm
 \label{FCP1}
\end{figure}

Figure \ref{FP22} shows in the enlarged scale first two FBP
profiles for this RP--$\frac 12$ that are rather close to the CP
one. The interfaces are
 $$
 y_{01} \sim 2.07 \quad \mbox{and} \quad y_{02} \sim 3.36.
  $$

\begin{figure}
\centering
\includegraphics[scale=0.75]{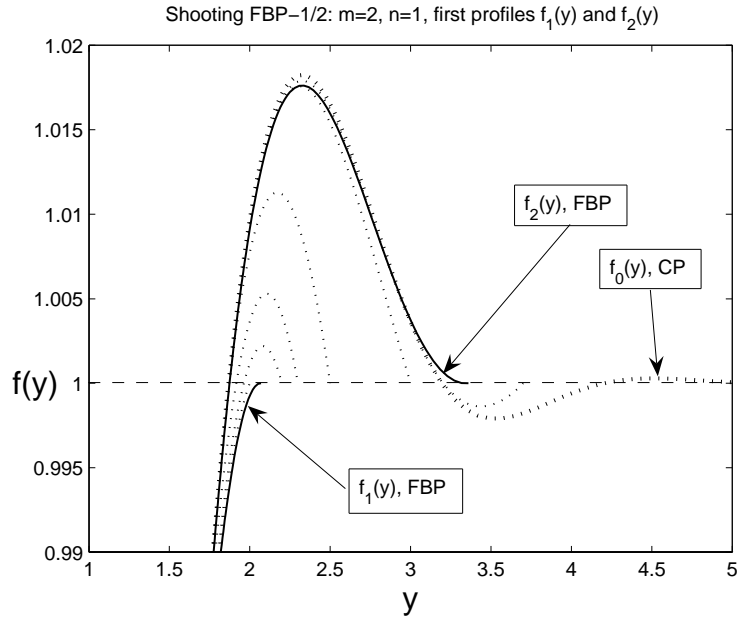}  
\vskip -.3cm
 \caption{\small Two first FBP solutions of the RP--$\frac 12$
(\ref{4}), (\ref{6}), (\ref{5R}) for $n=1$.}
   \vskip -.3cm
 \label{FP22}
\end{figure}

The key aspects of homotopic connections and branching at $n=0$
for the RP--$\frac 12$ remain the same as for the RP--1. The only
change is that,  convoluting with Heaviside data (\ref{2NN})
yields for $n=0$ the following similarity profile:
 \beq
 \label{ff11}
  \mbox{$
  f_0(y)= \int\limits_{-y}^{+\infty} F(z) \, {\mathrm d}z
   \quad (n=0).
 $}
   \eeq

\section{Interface equations for the
 TFE--4 with unstable linear diffusion}
 \label{Sect56}

We consider the TFE (\ref{5m}), and, studying finite propagation
close to $u=1$, we set
 $$
 1-u=v
 $$
  to obtain, up to small perturbations, the following unstable
  TFE:
   \beq
   \label{t1}
   v_t=- (|v| v_{xxx})_x - v_{xx}.
    \eeq
As usual, we take the simplest travelling wave (TW) solutions
 \beq
 \label{TW1}
 u(x,t)= f(y), \quad y= x- \l t \quad \Longrightarrow \quad -\l
 f'= -(|f| f''')'- f''.
  \eeq
  Hence, on integration once assuming the zero flux condition at the interface at
  $y_0=0$,
   \beq
   \label{t2}
   \l f= |f| f''' + f' \quad \mbox{for} \quad y>0, \quad f(0)=0.
    \eeq

\subsection{The FBP coincides with the CP}

In the first approximation, we neglect the non-stationary term on
the left-hand side of (\ref{t2}) and  consider the
 equation
  \beq
  \label{t4}
  |f| f''' + f' =0 \quad \mbox{for} \quad y>0, \quad f(0)=f'(0)=0.
   \eeq
If $f'(0)=0$, then the flux vanishes due to this ODE.  Taking the
expansion
 $f(y) = y^2 \rho(y)$
 and substituting into (\ref{t4}) we get
 $
 2 \rho'=-1 +... \, .
 $
 This gives the following first term expansion:
  \beq
  \label{t6}
  f(y) = y^2 | \ln y|+... \, .
 \eeq

Consider next the full equation (\ref{t2}). Dividing by $|f|$ and
integrating yields
 \beq
 \label{t7}
 \l y= f'' + \ln |f| + C \quad (C \in \re).
  \eeq

  \begin{proposition}
  \label{Pr1}
 ${\rm (i)}$ Equation $(\ref{t7})$ does not admit oscillatory solutions near the interface at
 $y=0$; and

 ${\rm (ii)}$ $(\ref{t7})$ does not have a decaying solution with $f(+\infty)=0$.
  \end{proposition}

 In fact, (ii) shows that (\ref{t7}) does not admit infinite propagation (for TWs of course).

 \smallskip

  \noi{\em Proof.} (i) Indeed, assuming that $f(y) \to 0$ as $y \to 0$, we have $\ln |f(y)| \ll -1$,
  so that $f''(y) >0$ and $f(y)$ is strictly convex in a neighbourhood.
  (ii) is straightforward.
  $\qed$

  \smallskip

  This means that the asymptotics (\ref{t6}) is already the actual
maximal regularity  behaviour corresponding to both the FBP and
CP, which thus coincide for such a PDE.

The dynamic equation for the interface at $x=s(t)$ follows from
(\ref{t7}) on differentiation, so we obtain the following equation
with a third-order interface operator $S$:
 \beq
 \label{t8}
 \mbox{$
 s'(t) = \l= (f'' + \ln |f|)'  =S[u] \equiv v_{xxx} + \frac {v_x}v
 \quad \mbox{at} \quad x=s(t).
  $}
  \eeq
  A rigorous justification of such equations assumes  using the von Mises variable
  $X(v,t)$ locally near sufficiently smooth and monotone interface,
  $$
  X(v(x,t),t) = x,
   $$
   though it leads to a number of technical difficulties; see \cite{Gl4}.

\section{Interface equations for the unstable nonlinear diffusion}
  \label{Sect5}

We now treat the interface propagation for the equation
(\ref{5m1N}), which, for $v=1-u \approx 0$, reduces to the
standard TFE--4 with the unstable diffusion
 \beq
 \label{h1}
 v_t= -(|v|v_{xxx})_x - (|v|^m v_x)_x.
  \eeq
  The TWs as in (\ref{TW1}), on integration, yield the ODE
 \beq
 \label{h2}
 \l f= |f| f''' + |f|^m f'.
  \eeq

\subsection{The FBP}

\smallskip

\noi \underline{\sc Case I.}
  The local analysis close to the interface at $y=y_0$ is pretty standard.
  Namely, we have got the pure TFE expansion like (\ref{k1})
  \cite{BPW, FB0},
 \beq
 \label{h3}
 \mbox{$
 f(y) = C y^2 + C_1 y^3 + ... \, ,
 \quad \mbox{where} \quad C_1= \frac \l 6, \quad C>0,
  $}
  \eeq
   provided that
   $$
    \mbox{$
   m > \frac 12.
    $}
   $$
This defines the standard interface equation
 \beq
 \label{h4}
 s'(t)= \l = 6 C_1 \equiv v_{xxx} \quad \mbox{at \,\, $x=s(t)$}.
  \eeq

\smallskip

 \noi \underline{\sc Case II.} For the critical exponent $m=
\frac 12$, all three terms in (\ref{h2}) are involved in the
expansion yielding in (\ref{h3})
 $$
  \mbox{$
C_1= \frac \l 6 - \frac 13 \, C^{m}.
 $}
 $$
 Therefore, the interface equation for TWs contains two operators
\beq
 \label{h5}
  \mbox{$
 s'(t)= \l =S[u] \equiv v_{xxx}+ 2 \bigl(\frac 12 \, v_{xx}\bigr)^{m} \quad \mbox{at \,\, $x=s(t)$}.
  $}
  \eeq

  \smallskip

\noi \underline{\sc Case III.} Finally, for $m  \in (0, \frac
12)$, the FBP expansion is not $C^3$-smooth. For instance, in the
first such range $m \in (\frac 14,\frac 12)$, the expansion is
 \beq
 \label{h6}
  \mbox{$
 f(y) = C y^2 + C_1 y^{2+2m} + \frac \l 6 \, y^3 +  ... \, ,
 \quad \mbox{where} \quad
  C_1=-  \frac 1{2m(1+m)(1+2m)} \, C^m,
   $}
  \eeq
  and, in the second term, the exponent satisfies
   $
   2<2+2m <3.
    $
Therefore, $f \not \in C^3([0,1])$, so the equality
 $\l=f'''(0)$ makes no sense and the interface
equation becomes rather tricky. The simplest way to resolve
(\ref{h6}) with respect to the speed $\l=s'(t)$ is to decompose
the TW solutions as follows:
 $$
 v(x,t)= \hat v(x,t) + \tilde v(x,t),
 $$
 where $\hat v(x,t)$ is $C^3$-smooth at the interface and $\tilde
 v(x,t)$ is not. Then the interface equation in this non-regular
 case takes the standard (similar to (\ref{h4})) form
  $$
  s'(t)=(\hat v)_{xxx} \quad \mbox{at} \,\,\,\, x= s(t).
  $$
 An alternative, higher-order version can be obtained by observing
 from (\ref{h2}) that, for solutions (\ref{h6}),
  $$
   \mbox{$
  \l C y^2 + ... = f f''' + f^m f' \quad \bigl(C= \frac 12 \,
  f''(0)\bigr).
   $}
 $$
Differentiating two times yields the following interface equation:
 $$
  \mbox{$
  s'(t)=S[u] \equiv \frac 1{v_{xx}} \,(v v_{xxx}+ v^m v_x)_{xx} \quad \mbox{at} \,\,\,\, x= s(t).
   $}
  $$
The case $m=\frac 14$ is critical, where the second term on the
right-hand side
 in (\ref{h2}) adds an extra operator into the interface equation
 that now reads (the above differential form is also available)
 $$
  \mbox{$
  s'(t)=(\hat v)_{xxx}- \frac 83 \bigl(\frac 12 v_{xx}\bigr)^{2m-1} \quad \mbox{at} \,\,\,\, x= s(t).
   $}
  $$

A rigorous proof of such types of interface equations is a
difficult open problem.


\subsection{The CP}

\smallskip

\noi \underline{\sc Case I.} For $m > \frac 13$, the propagating
($\l<0$) TW solutions are oscillatory near interfaces and have the
behaviour as in (\ref{k3}), i.e.,
 \beq
 \label{h7}
 f(y) = y^3 \varphi(s), \quad s= \ln y,
  \eeq
  where $\varphi(s)$ is a periodic solution of the ODE (\ref{k4})
   with $\l_0 \mapsto -\l$.

   \smallskip

\noi \underline{\sc Case II.} In the critical case $m= \frac 13$,
we still have the expansion (\ref{h7}), where the oscillatory
component solves a different ODE,
 \beq
 \label{h8}
 P_3(\varphi)= \l \, {\rm sign} \, \varphi - (\varphi'+3
 \varphi)|\varphi|^{-\frac 23}.
  \eeq
Existence-uniqueness of a stable periodic orbit becomes harder and
remains an open problem.  Figure \ref{F56}(a) shows the unique
stable periodic solution of (\ref{h8}). In Figure (b), for
comparison, we show the periodic behaviour in the ODE for the
stable counterpart of the TFE with $m= \frac 13$,
 \beq
 \label{h9}
 v_t = - (|v| v_{xxx})_x + (|v|^{\frac 13} v_x)_x
  \,\, \Longrightarrow
  \,\, P_3(\varphi)= \l \, {\rm sign} \, \varphi + (\varphi'+3
 \varphi)|\varphi|^{-\frac 23}.
 \eeq

\begin{figure}
\centering
\subfigure[unstable diffusion, $m=\frac 13$]{
\includegraphics[scale=0.52]{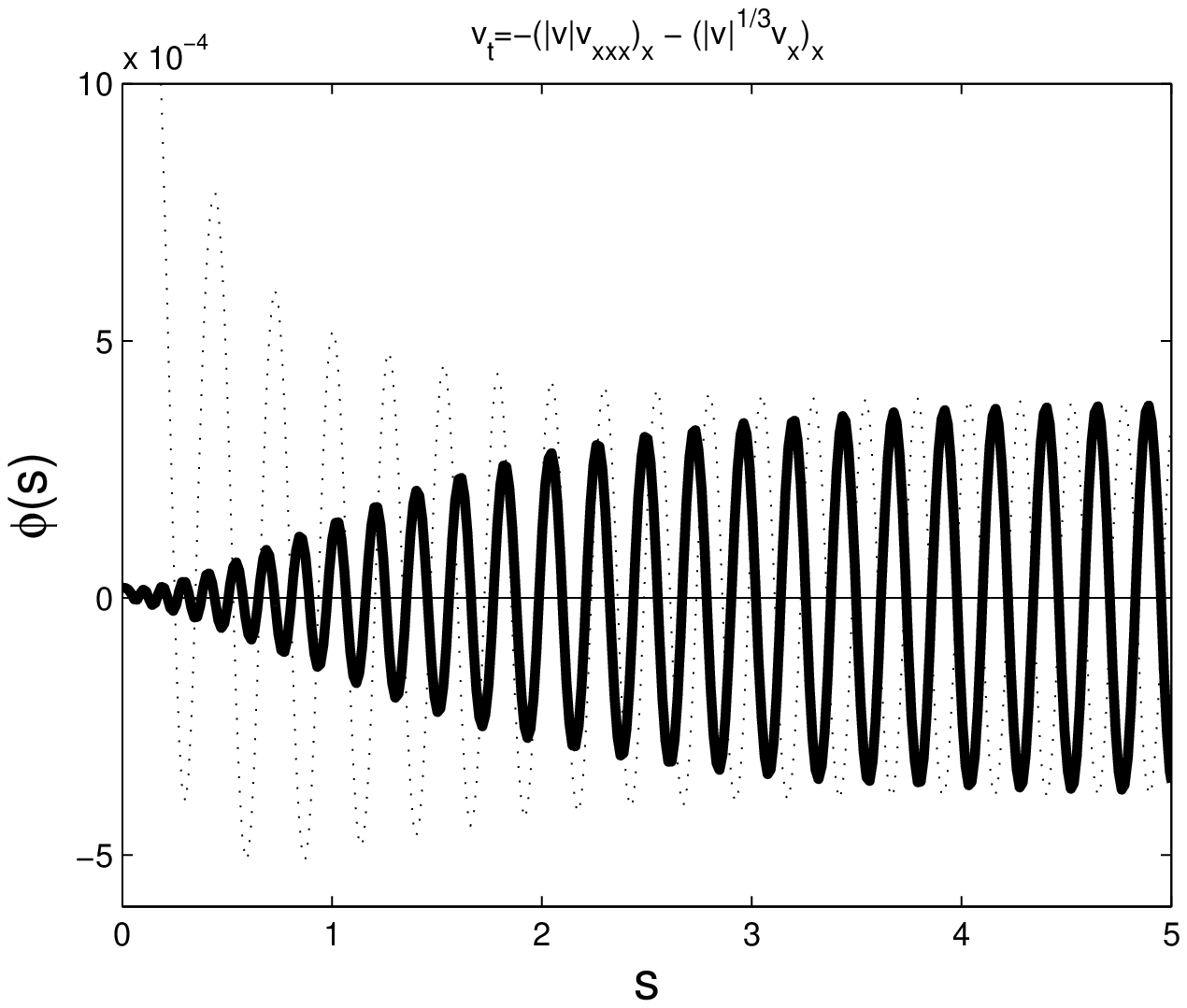} 
}
\subfigure[stable diffusion, $m=\frac 13$]{
\includegraphics[scale=0.52]{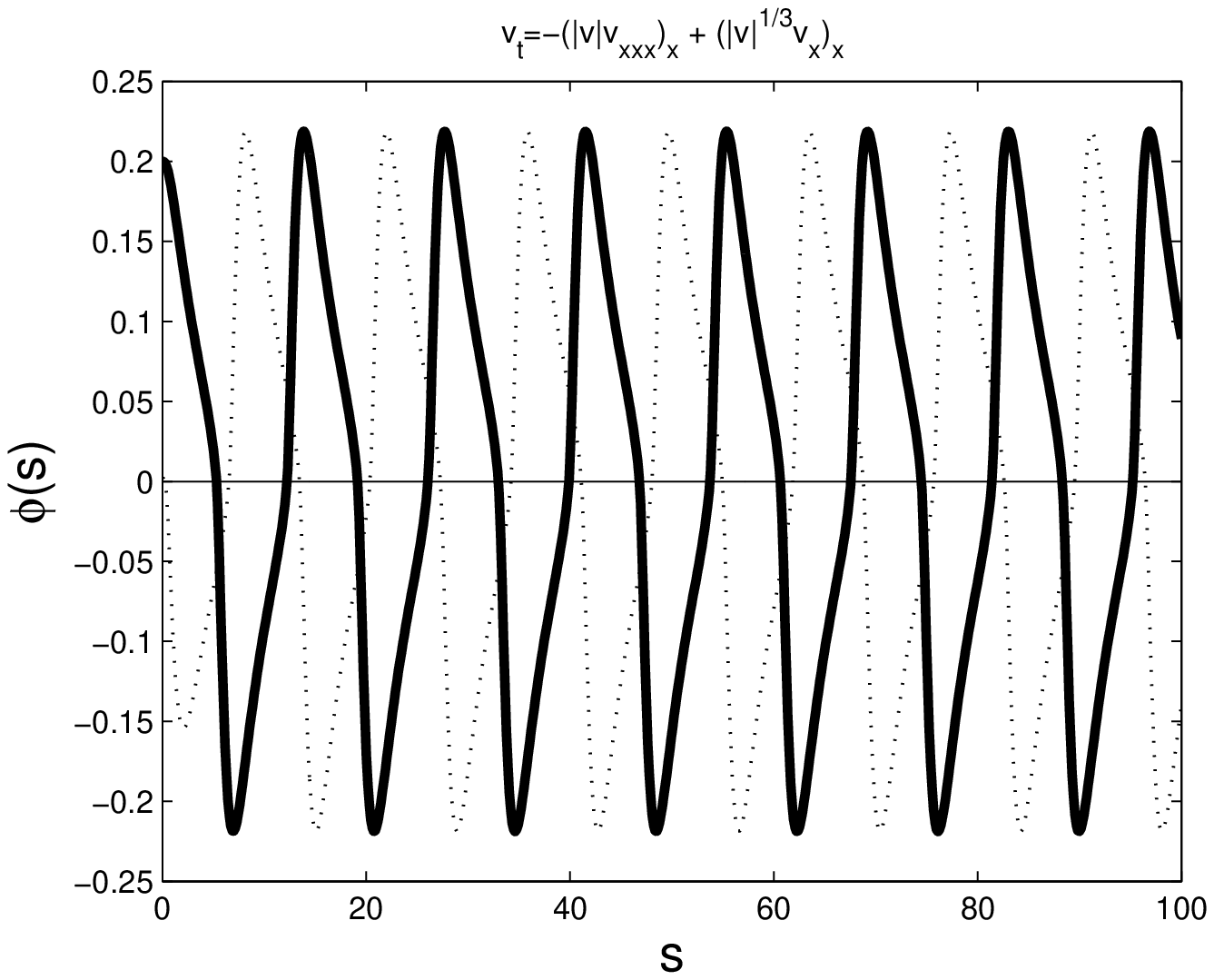} 
}
   \vskip -.3cm
    \caption{\small  Stable periodic changing sign behaviour of the oscillatory component
    for the ODE (\ref{h8}) (a), and (\ref{h9}) (b); $\l =-1$.}
   \vskip -.3cm
 \label{F56}
\end{figure}

\smallskip

\noi \underline{\sc Case III.} Finally, in the range
 $$
 \mbox{$
 0<m< {\rm min} \,\big\{\frac 13,n\big\}
 $}
 $$
two terms on the right-hand side of (\ref{h2}) are leading that
yields the same quadratic expansion (\ref{h6}) as in the FBP. It
seems that,  in this parameter range (and also for $m=0$ in
Section \ref{Sect56}), the CP and FBP settings coincide.


In the CP, in view of the oscillatory nature of solutions in most
cases, even formal derivation of interface equations becomes
difficult, to say nothing about their  rigorous justification.





\section{On higher-order TFEs: the RP--1 for $m=3$}
 \label{SectHi}

Similarity solutions persist for higher-order TFEs (\ref{4m}). The
matching becomes more-parametric and the topology of local bundles
gets rather delicate. Numerics also get essentially more involved.

\subsection{On the CP profile and oscillations}

As an example, in Figure \ref{Fm3}, we present the CP profile
satisfying (\ref{4m}) for the case $m=3$ and $n=1$, i.e., for the
{\em sixth-order TFE} (TFE--6). Instead of (\ref{k3}), the
oscillatory behaviour about $f(y) \equiv 1$ is smoother,
 \beq
 \label{k33}
 f(y)=1- (y_0-y)^5 \varphi(s)+... \, , \quad s= \ln(y_0-y),
  \eeq
where $\varphi(s)$ solves a fifth-order ODE
 \beq
 \label{kk1}
  P_5(\varphi) \equiv \varphi^{(5)}+15\varphi^{(4)}
+85\varphi''' 
 + 255\varphi''
+\, 274\varphi' + 120\varphi = \l_0 \, {\rm sign} \, \varphi,
 \eeq
 where $\l_0= \frac 1{6}\, y_0$.
See details on derivation  in \cite[\S~12.2]{GBl6}. Figure
\ref{Figm3} shows the unstable periodic solution of (\ref{kk1})
for $\l_0=1$, which by (\ref{k33})  describes the oscillatory
behaviour near the interface at
 \beq
 \label{y0CP1}
 y_{0,{\rm CP}} \sim 13.
  \eeq



\begin{figure}
\centering
\subfigure[CP profile]{
\includegraphics[scale=0.52]{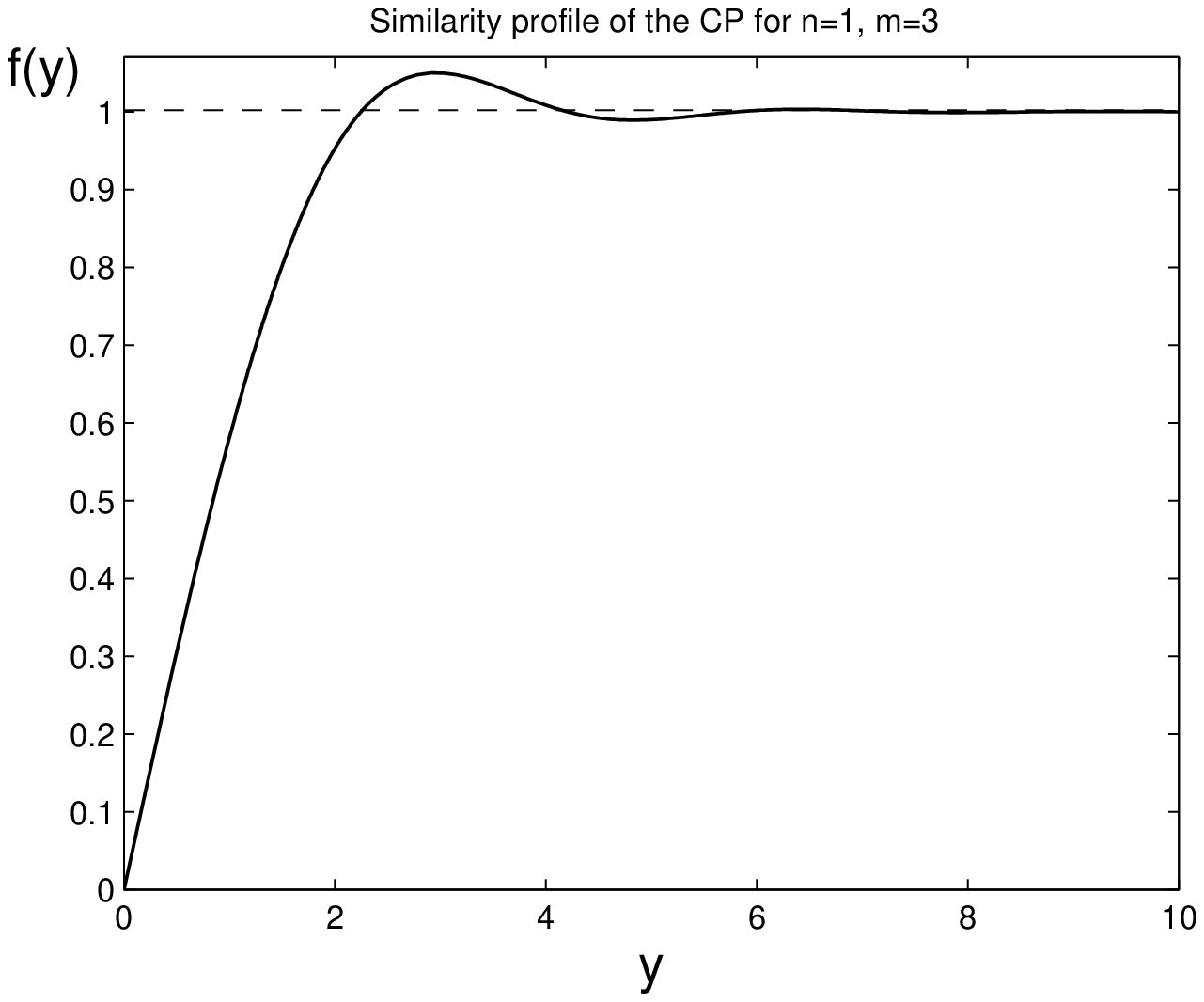}
}
\subfigure[oscillations, enlarged]{
\includegraphics[scale=0.52]{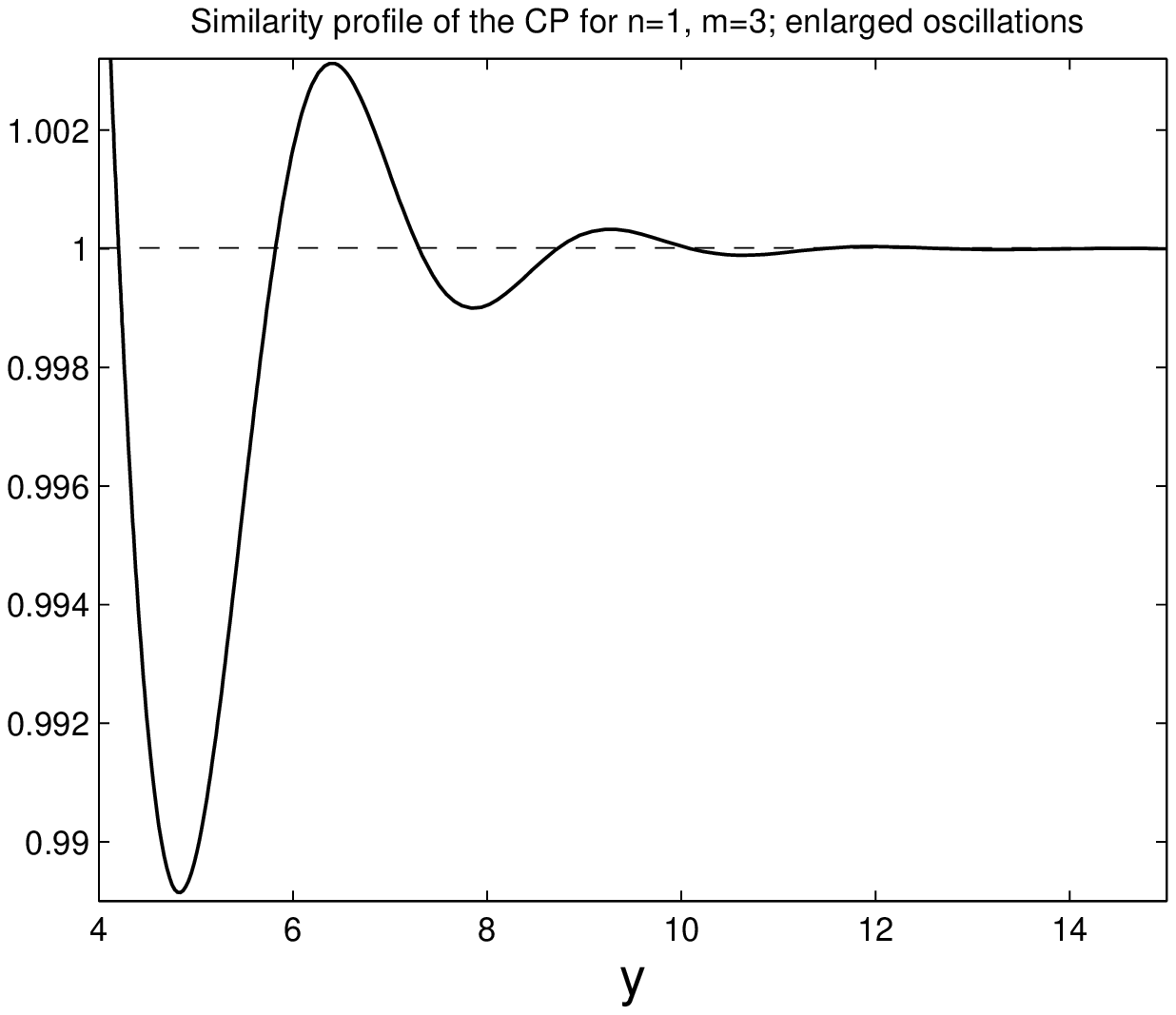}
}
   \vskip -.3cm
    \caption{\small  The unique CP profile $f(y)$
    solving (\ref{4m}), (\ref{7}) for $m=3$, $n=1$; the profile (a), and
    enlarged oscillations about 1 (b).}
   \vskip -.3cm
 \label{Fm3}
\end{figure}

\begin{figure}
\centering
\includegraphics[scale=0.60]{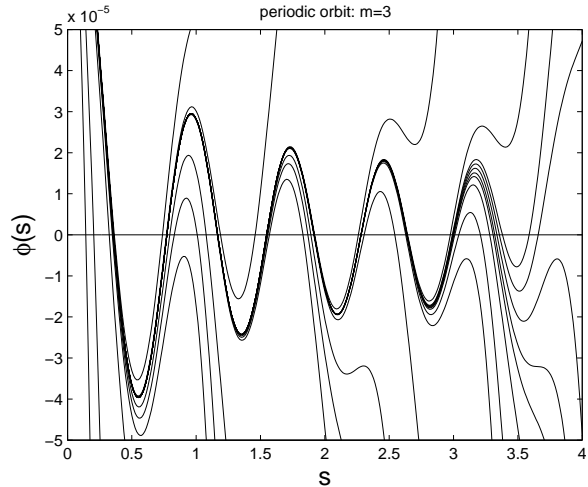}  
\vskip -.3cm \caption{\small The trace of the unstable periodic
orbit of (\ref{kk1}), with $\l_0=1$ ($y_0=6$).}
   \vskip -.3cm
 \label{Figm3}
\end{figure}

 \subsection{Shooting FBP profiles}

The results of a standard shooting of the first FBP profiles are
presented in Figure \ref{Figm39}. Similar to shooting in Figure
\ref{F01}, here we have chosen the same strategy  for higher-order
equations: for $m=3$, we fix three free-boundary conditions from
(\ref{5m}) excluding the second one, so changing the interface
location $y_0$ as the free parameter, we want to satisfy the zero
contact angle condition (\ref{77}).
As a consequence, we see in Figure \ref{Figm39} a consistent
illustration of the limit FBP--CP process indicated earlier in the
schematic Figure \ref{Fcount}. Namely, we obtain four first FBP
profiles, $f_1(y)$, $f_2(y)$, $f_3(y)$, and $f_4(y)$, with
interfaces positions at
 \beq
 \label{int11}
 y_{01} \sim 3.05, \,\,\,
  y_{02} \sim 4.85, \,\,\,
   y_{03} \sim 6.55, \,\,\, \mbox{and} \,\,\,
    y_{04} \sim 8.1.
    \eeq
These are smaller than the Cauchy one (\ref{y0CP1}), which is
expected to be the limit of all other FBP interfaces $\{y_{0k}\}$
as $k \to \infty$.

As usual, only the first profile  $f_1(y)$ satisfies the physical
range condition (\ref{range1}). Other profiles are oscillatory
about the constant equilibrium $f=1$, and as it is seen in Figure
\ref{Figm39}, rather thoroughly mimic the oscillations of the
linear solution (the boldfaced dotted line) defined by the same
formula (\ref{b3}).
 This figure
confirms that Conjecture \ref{Sect31}.2 applies also to the TFE--6
(and seems to any of $2m$th-order ones (\ref{8})).

\begin{figure}
\centering
\includegraphics[scale=0.8]{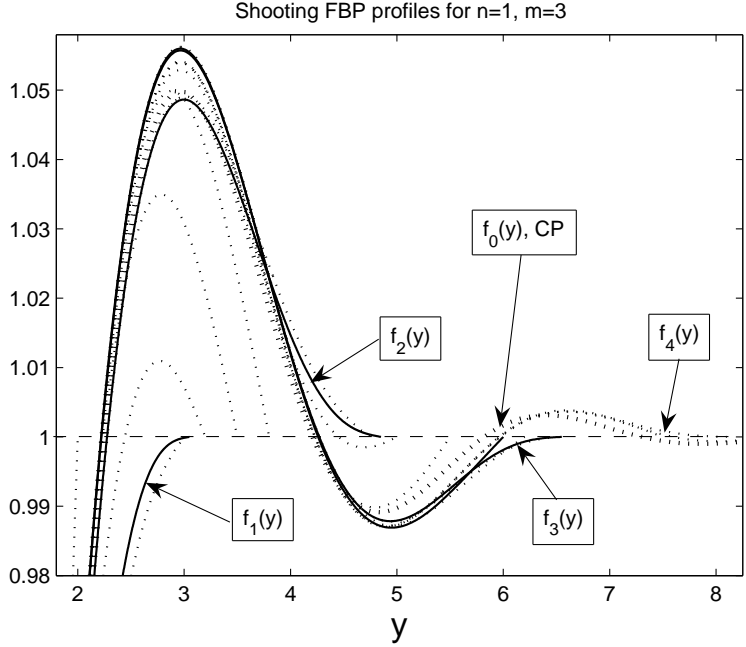}
\vskip -.3cm
 \caption{\small Shooting  first four FBP profiles
  satisfying (\ref{4m}), (\ref{5m}) for $m=3$,
$n=1$.}
   \vskip -.3cm
 \label{Figm39}
\end{figure}

For comparison, in Figure \ref{Fig00}, we present first three FBP
profiles for the linear case $n=0$, i.e., $f$ solves the FBP
 \beq
 \label{ff1}
 \begin{matrix}
 f^{(6)} + \frac 16 \, y f'=0 \quad \mbox{for}
 \quad y>0,\qquad \smallskip\smallskip \\
 f(0)=f''(0)=f^{(4)}(0)=0, \quad f(y_0)=1, \,\,\,
 f'(y_0)=f''(y_0)=f^{(5)}(y_0)=0.\qquad
  \end{matrix}
  \eeq
  The interface positions are larger than that in (\ref{int11}) for
  $n=1$,
 $$
 y_{01} \sim 3.29, \,\,\,
  y_{02} \sim 6.57, \,\,\,\mbox{and} \,\,\,
   y_{03} \sim 9.72.
    $$
Observe that the FBP profiles are close to corresponding humps of
the CP profile denoted by the boldface dotted line. The general
geometry of the mutual location of FBP and the CP profiles for
$n=0$ is quite similar to the ``nonlinear" one in Figure
\ref{Figm39} for $n=1$ (the  figures are topologically
equivalent).

\begin{figure}
\centering
\includegraphics[scale=0.8]{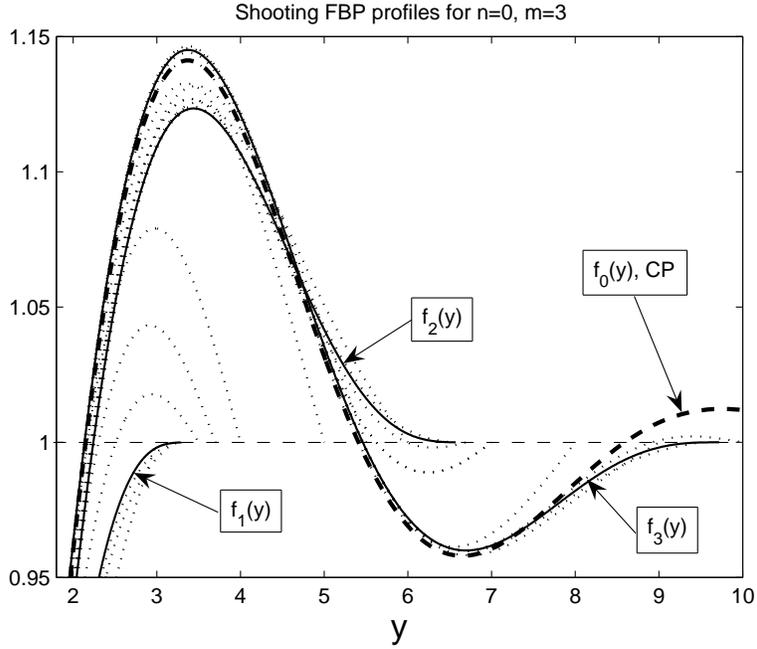}
\vskip -.3cm
 \caption{\small Shooting  first three FBP profiles
  satisfying the linear FBP  problem (\ref{ff1}).}
   \vskip -.3cm
 \label{Fig00}
\end{figure}

\subsection{Branching at $n=0$}

For $n=0$, (\ref{8}) for $m=3$ becomes the {\em tri-harmonic
equation}
 \beq
 \label{l1S}
 u_t=u_{xxxxxx} \quad \mbox{in} \quad \re \times \re_+.
  \eeq
Its {fundamental solution}  is given by
 \beq
 \label{b1S}
  \mbox{$
 b(x,t) = t^{-\frac 16} F(y), \quad y =  \frac x{t^{ 1/6}},
  $}
  \eeq
  where $F$ is the unique symmetric solution of the problem
   \beq
   \label{b2S}
   \textstyle{
 {\bf B} F \equiv   F^{(6)}+ \frac 16 \, F' y + \frac 16\, F=0 \quad \mbox{in}
   \quad \re, \quad \int F=1.
   }
 \eeq
Extension of the similarity profiles of the ODE (\ref{4m}) from
the branching point $n=0$ is the same as that for the TFE--4 in
Section \ref{SectExi}. Spectral properties of the corresponding
linear operators ${\bf B}$ in (\ref{b2S}) and the adjoint one
 $$
 \mbox{$
  {\bf B}^*=
 D_y^{6} - \frac 1{6}\, y D_y
 $}
 $$
 in $L^2_\rho(\re)$ and $L^2_{\rho^*}(\re)$ respectively
 can be found in \cite{Eg4}.

  In Figure \ref{FigBr4},
we show branching of the similarity profiles at $n=0$ (the  dotted
line), where, with the step $\D p=0.2$, we cover the range $p \in
[-0.8,1]$. All the profiles are oscillatory, but for $n<0$  the
interface is situated at $y_0=+\infty$.

\begin{figure}
\centering
\includegraphics[scale=0.70]{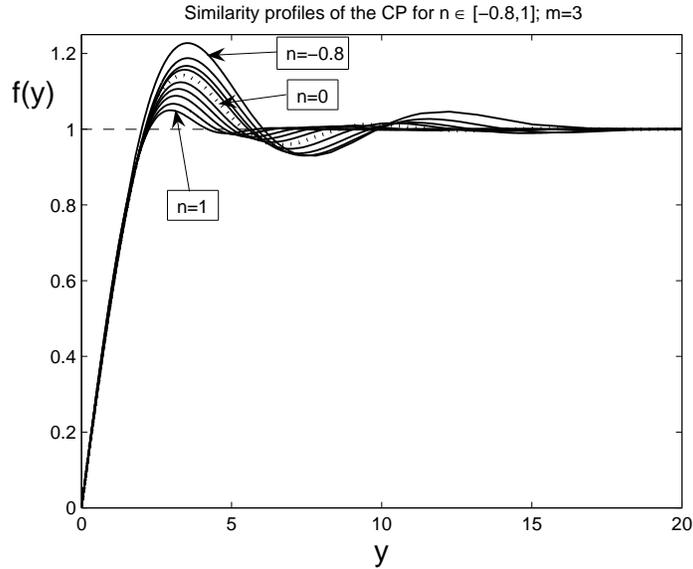}
\vskip -.3cm
 \caption{\small Similarity profiles of the CP satisfying (\ref{4m}) for $m=3$
 appearing at  the branching point $n=0$ (the dotted line).}
   \vskip -.3cm
 \label{FigBr4}
\end{figure}

\smallskip

{\bf Acknowledgements.} The author would like to thank
A.~Novick-Cohen
 and A.~Shishkov for  discussions of the physical motivation of the problem for (\ref{eq1})
 and of some mathematical results respectively.


\end{document}